\documentclass[10pt]{amsart}

\usepackage[mathscr]{eucal}
\usepackage{amsthm,amscd,latexsym,xspace,amssymb}

\newcommand{\Supp}{{\rm{Supp}\ts}}
\newcommand{\non}{\nonumber}
\newcommand{\wt}{\widetilde}

\newcommand{\ot}{\otimes}
\newcommand{\la}{\lambda}
\newcommand{\al}{\alpha}
\newcommand{\be}{\beta}
\newcommand{\ga}{\gamma}
\newcommand{\Ga}{\Gamma}

\newcommand{\ts}{\,}
\newcommand{\tss}{\hspace{1pt}}
\newcommand{\U}{ {\rm U}}
\newcommand{\Y}{ {\rm Y}}

\newcommand{\Z}{\mathbb{Z}\tss}
\newcommand{\N}{\mathbb{N}\tss}

\newcommand{\gl}{\mathfrak{gl}}

\newcommand{\gr}{ {\rm gr}}

\newcommand{\dmo}{\text{-mod}}
\newcommand{\bra}{{\rm (\ts}}
\newcommand{\ket}{{\rm ) }}
\newcommand{\mP}{\mathcal P}

\newtheorem{thm}{Theorem}[section]
\newtheorem{prop}[thm]{Proposition}
\newtheorem{cor}[thm]{Corollary}
\newtheorem{lem}[thm]{Lemma}

\theoremstyle{definition}
\newtheorem{defin}[thm]{Definition}

\theoremstyle{remark}
\newtheorem{remark}[thm]{Remark}
\newtheorem{example}[thm]{Example}

\newcommand{\bth}{\begin{thm}}
\renewcommand{\eth}{\end{thm}}
\newcommand{\bpr}{\begin{prop}}
\newcommand{\epr}{\end{prop}}
\newcommand{\ble}{\begin{lem}}
\newcommand{\ele}{\end{lem}}
\newcommand{\bco}{\begin{cor}}
\newcommand{\eco}{\end{cor}}
\newcommand{\bde}{\begin{defin}}
\newcommand{\ede}{\end{defin}}
\newcommand{\bex}{\begin{example}}
\newcommand{\eex}{\end{example}}
\newcommand{\bre}{\begin{remark}}
\newcommand{\ere}{\end{remark}}
\newcommand{\bal}{\begin{aligned}}
\newcommand{\eal}{\end{aligned}}
\newcommand{\beq}{\begin{equation}}
\newcommand{\ben}{\begin{equation*}}

\def\beql#1{\begin{equation}\label{#1}}

\newcommand\dcl\DeclareMathOperator

\dcl\Sp{Specm}

\newcommand\1{{\mathbf 1}}
\dcl\cfs{cfs}
 \dcl\supp{supp}
\newcommand\xar[1]{\xrightarrow{\quad #1 \quad }}
\dcl\Ker{Ker}
\dcl\Hom{Hom}
\dcl\Ext{Ext}
\dcl\Ann{Ann}
\dcl\Ob{Ob}
\dcl\im{Im}
\dcl\mo{mod}
\dcl\rank{rank}
\dcl\M{M}

\renewcommand{\L}{\mathrm {L}}
\newcommand{\G}{\mathrm {G}}
\newcommand{\mL}{\mathcal L}
\newcommand{\bL}{\mathbb L}
\newcommand{\bW}{\mathbb W}
\newcommand{\bm}{\mathbf m}
\newcommand{\bn}{\mathbf n}
\renewcommand{\o}{\overline{0}}
\renewcommand{\t}{\overline{t}}
\newcommand\T{\mathrm {T}}

\begin{document}

\title{HARISH-CHANDRA MODULES FOR  YANGIANS}
\author{Vyacheslav Futorny}
\address{Institute of Mathematics and Statistics\\
University of S\~ao Paulo,\\
Caixa Postal 66281- CEP 05315-970\\
S\~ao Paulo, Brazil}
\email{futorny@ime.usp.br}
\author{Alexander Molev}
\address{
School of Mathematics and Statistics\\
University of Sydney,
NSW 2006, Australia}
\email{alexm@maths.usyd.edu.au}
\author{Serge Ovsienko}
\address{
Faculty of Mechanics and Mathematics\\
Kiev Taras Shevchenko University\\
Vla\-di\-mir\-skaya 64, 00133, Kiev, Ukraine}
\email{ovsienko@sita.kiev.ua }

\begin{abstract}
We study Harish-Chandra  representations of the Yangian
$\Y(\gl_2)$ with respect to a natural
maximal commutative subalgebra which satisfy a polynomial
condition. We  prove an analogue of the Kostant theorem showing that
the restricted Yangian $\Y_p(\gl_2)$ is a free module over
the corresponding subalgebra
$\Gamma$ and show that every character of $\Gamma$ defines a
finite number of irreducible Harish-Chandra modules over $\Y_p(\gl_2)$.
We  study the
categories of generic Harish-Chandra modules, describe their
simple  modules and indecomposable modules in tame blocks.
\end{abstract}

\vspace{0,2cm}

\maketitle


%%%%%%%%% provisorisch

\tableofcontents

\section {Introduction}
\setcounter{equation}{0}

Throughout the paper we fix an algebraically closed
field {${\Bbbk}$}{} of characteristic $0$.
Consider the pair $(U,\Gamma)$ where
$U$ is an associative $\Bbbk$-algebra and
$\Gamma$ is a subalgebra of $U$. Denote by $\cfs\Gamma $ the
{\em cofinite spectrum\/}
of $\Gamma$, i.e.,
\ben
\cfs \Gamma =\{\text{maximal two-sided ideals $\bm$
of $\Gamma$}\ |\ \dim \Gamma/\bm<\infty\}.
\end{equation*}
A finitely generated module $M$ over $U$ is called
a {\em Harish-Chandra module\/} (with respect to $\Gamma$)
if
\ben
M=\underset{{\bm} \in \cfs {\Gamma}}{\bigoplus}M({\bm}),
\end{equation*}
where
\ben
M({\bm}) \ = \ \{ x\in M\ | \ {\bm}^k x =0\quad
\text{for some}\quad k\geq 0\}.
\end{equation*}

Harish-Chandra  modules play a central role
in the classical representation
theory; see e.g. Dixmier~\cite{d:ae}. In particular, weight modules
over a semisimple Lie algebra
are Harish-Chandra modules with respect to a Cartan subalgebra.
Another important example is provided by
the Gelfand--Tsetlin modules \cite{dfo:gz} over the
universal enveloping algebra $U(\gl_n)$ of the general linear Lie
algebra $\gl_n$. They are Harish-Chandra modules with respect to the
Gelfand--Tsetlin subalgebra of $U(\gl_n)$. The latter is the
commutative subalgebra generated by the centers of $U(\gl_k)$,
$k=1, \ldots,n$.  A theory of Harish-Chandra modules
for general pairs $(U,\Gamma)$ is developed in \cite{dfo:hc}.

An irreducible Harish-Chandra module $M$ is said to be {\em extended\/}
from $\bm\in \cfs \Gamma$ if $M(\bm)\neq 0$.
A central problem in the theory of Harish-Chandra
modules is to investigate the existence and uniqueness
conditions for such an extension.
In the case where the extension is unique, the irreducible
Harish-Chandra modules are parametrized by some equivalence classes
of the elements of
$\cfs \Gamma$. It has been recently proved in \cite{o:fs} that
for the case of Gelfand--Tsetlin modules over
$\gl_n$  the number of pairwise
non-isomorphic irreducible modules extended from a given $\bm\in
\cfs\Gamma$ is nonzero and finite.

In this paper we begin a detailed study of Harish-Chandra
modules over the Yangians.
The {\it Yangian for\/} $\gl_n$ is a unital associative algebra
$\Y(\gl_n)$ over $\Bbbk$ with countably many generators
$t_{ij}^{(1)},\ t_{ij}^{(2)},\dots$ where $1\leq i,j\leq n$, and
the defining relations
\beql{defrel} (u-v)\ts
[t_{ij}(u),t_{kl}(v)]=t_{kj}(u)\ts t_{il}(v)-t_{kj}(v)\ts
t_{il}(u),
\end{equation}
where
\ben
t_{ij}(u) = \delta_{ij} + t^{(1)}_{ij} u^{-1} + t^{(2)}_{ij}u^{-2} +
\cdots
\end{equation*}
and $u,v$ are formal variables. This algebra originally appeared
in the works on the {\it quantum inverse
scattering method\/}; see e.g. Takhtajan--Faddeev~\cite{tf:qi},
Kulish--Sklyanin~\cite{ks:qs}. The term ``Yangian" and generalizations
of  $\Y(\gl_n)$ to an arbitrary simple Lie algebra were introduced by
Drinfeld~\cite{d:ha}. He then classified finite-dimensional irreducible modules
over the Yangians in \cite{d:nr} using earlier results of Tarasov~\cite{t:sq, t:im}.
An explicit construction of every such
module over $\Y(\gl_2)$ is given in those papers by Tarasov and also
in the work by Chari and Pressley~\cite{cp:yr}.
Apart from this case, the structure
of a general finite-dimensional irreducible representation of the Yangian
remains unknown.
In the case of $\Y(\gl_n)$ a description of
{\em generic\/} modules was given in \cite{m:gt} via Gelfand--Tsetlin bases.
A more general class of {\em tame\/} representations of $\Y(\gl_n)$
was introduced and explicitly constructed by Nazarov and Tarasov~\cite{nt:ry}.
Another family of representations has been
constructed in \cite{m:ic} via tensor products
of the so-called {\em evaluation modules\/}.
An important role in these works is played by the
{\it Drinfeld generators\/}~\cite{d:nr}
\ben
a_i(u),\quad i=1,\dots,n,\qquad b_i(u),\ \  c_i(u),\quad i=1,\dots,n-1
\end{equation*}
of the algebra $\Y(\gl_n)$ which are defined as certain {\it
quantum minors\/} of the matrix $T(u)=\big(t_{ij}(u)\big)$. The
coefficients of the series $a_i(u)$, $i=1,\dots,n$ form a
commutative subalgebra of $\Y(\gl_n)$ which can be regarded as an
analogue of the Gelfand--Tsetlin subalgebra of
$\U(\gl_n)$.
We shall consider the Harish-Chandra modules for $\Y(\gl_n)$
with respect to this particular subalgebra.
So, the Harish-Chandra modules for $\Y(\gl_n)$ are natural
analogues of the Gelfand--Tsetlin modules for $\gl_n$ \cite{dfo:gz}.
Note also that the tame modules over $\Y(\gl_n)$ \cite{nt:ry}
is a particular case of
Harish-Chandra modules.

In this paper we are concerned with Harish-Chandra modules for
the Yangian $\Y(\gl_2)$. Recall that every irreducible finite-dimensional
$\Y(\gl_2)$-module contains a unique vector $\xi$ annihilated by
$t_{12}(u)$ and which is an eigenvector for the Drinfeld generators
$a_1(u)$ and $a_2(u)$ defined by
\beql{a1u}
a_1(u)=t_{11}(u)\ts t_{22}(u-1)-t_{21}(u)\ts t_{12}(u-1),\qquad a_2(u)=t_{22}(u);
\end{equation}
see \cite{t:sq, t:im} and \cite{cp:yr}.
Moreover, there exists an automorphism
 $t_{ij}(u)\mapsto c(u)\ts t_{ij}(u)$ of $\Y(\gl_2)$,
where $c(u)\in 1 + u^{-1}\ts\Bbbk[[u^{-1}]]$,
such that the eigenvalues of $\xi$ become polynomials in $u^{-1}$
under the corresponding twisted action of the Yangian.
This prompts the introduction of the class of
{\it Harish-Chandra polynomial\/} modules  over  $\Y(\gl_2)$, i.e., such
 Harish-Chandra
modules where the operators $a_1(u)$ and $a_2(u)$ are polynomials.
More precisely, due to \eqref{a1u},
it is natural to require that for some positive integer $p$
the polynomials $a_1(u)$ and $a_2(u)$ have degrees $2p$ and $p$, respectively.
Note that $a_1(u)$ is the {\it quantum determinant\/} of the matrix $T(u)$ \cite{ik:lm},
\cite{ks:qs}.
Its coefficients are algebraically independent generators
of the center of $\Y(\gl_2)$.

We can interpret the definition of Harish-Chandra polynomial modules
using the algebra $\Y_p(\gl_2)$ called the {\it Yangian of level\/} $p$;
see Cherednik~\cite{c:ni, c:qg}. It is defined as the quotient
of $\Y(\gl_2)$ by the ideal generated by the elements  $t_{ij}^{(r)}$
with $r\geq p+1$. A Harish-Chandra polynomial module over $\Y(\gl_2)$
is just a Harish-Chandra module over $\Y_p(\gl_2)$ for some positive integer $p$.
In what follows we shall consider
Harish-Chandra modules over
$\Y_p(\gl_2)$ with respect to the commutative subalgebra $\Gamma$
generated by the coefficients of $a_1(u)$ and
$a_2(u)$.

Let us now describe our main results.
First, we prove that $\Y_p(\gl_2)$ is free as a left (right)
$\Gamma$-module (Theorem~\ref{thm:freedom}).
This is an analogue of the well-known Kostant theorem \cite{k:lg}.
Each character
of $\Gamma$ can therefore be extended
to an irreducible $\Y_p(\gl_2)$-module.
An important role in
our study is played by certain universal
Harish-Chandra modules over $\Y_p(\gl_2)$
(Theorem~\ref{thm:basis}) such
that every irreducible module in  {$ \mathbb H(\Y_p(\gl_2),{\Gamma})$}{}
is a quotient of the corresponding universal module.

Further, we show that $\Gamma$ is a Harish-Chandra subalgebra
(Theorem~\ref{thm:GT-HC}) in the sense of \cite{dfo:hc}
which allows us to apply the general theory of \cite{dfo:hc}
to the study of Harish-Chandra modules for $\Y_p(\gl_2)$.
In particular,
it provides an equivalence between the category
{$ \mathbb H(\Y_p(\gl_2),{\Gamma})$} of Harish-Chandra modules and
the category of finitely generated modules over a certain category
${\mathscr A}$ whose objects are the maximal ideals of $\Gamma$.
We then use this to prove that the number
of pairwise non-isomorphic extensions of a character of $\Gamma$
to an irreducible $\Y_p(\gl_2)$-module
is finite (Theorem \ref{thm:main}).
The full subcategory
{$ \mathbb HW(\Y_p(\gl_2),{\Gamma})$} of {$ \mathbb H(\Y_p(\gl_2),{\Gamma})$}
which consists
of modules with diagonalizable action of $\Gamma$
turns out to be equivalent to the category of
finitely generated modules over a certain quotient category of
${\mathscr A}$ (Section~\ref{subsec:Harish}).  In Section~\ref{sec:Generic}
we study a full
subcategory in  {$ \mathbb HW(\Y_p(\gl_2),{\Gamma})$}{} of generic
modules, this imposes a certain condition on the
eigenvalues of $a_2(u)$ while those of $a_1(u)$ are arbitrary. In
particular, we give a complete description of irreducible modules
(Theorem \ref{thm:suppar}) and indecomposable modules in tame
blocks of this category (Theorem \ref{thm:reptype}).

\newpage

\section{Preliminaries}
\setcounter{equation}{0}

\subsection{Harish-Chandra subalgebras}\label{subsec:Harish}
\label{subsec:hca}

In this paper we shall only consider the pairs $(U,\Gamma)$
where the subalgebra $\Gamma$ of $U$ is commutative.
In this case $\cfs \Gamma$
coincides with the set $\Sp \Gamma$ of all maximal ideals in
$\Gamma$.

We let $U\dmo$ denote the category of finitely generated
left modules over an associative algebra $U$.
The Harish-Chandra modules for the pair $(U,\Gamma)$ form a
full abelian subcategory in $U\dmo$ which we
denote by {$ \mathbb H(U,{\Gamma})$}{}.
A Harish-Chandra module $M$ is called {\em weight\/} if
the following condition holds:
for all ${\mathbf m}\in \Sp {\Gamma}$ and all $x\in
M({\mathbf m})$ one has ${\mathbf m}\ts x=0$.
The full subcategory of
{$ \mathbb H(U,{\Gamma})$}{} consisting of weight modules will be
denoted  {$ \mathbb HW(U,{\Gamma})$}{}. The {\em support\/} of a
Harish-Chandra module $M$ is the subset $\Supp M\subseteq$ $\Sp {\Gamma}$
which consists of those ${\mathbf m}$ which have the property
$M({\mathbf m})\ne 0$.
 If for a given $\mathbf m$ there exists an irreducible
Harish-Chandra module $M$ with $M({\mathbf m})\ne 0$ then we say
that $\mathbf m$ {\em extends\/} to $M$.

%The notion of a Harish-Chandra subalgebra \cite{dfo:hc} is an
%effective tool for the study of the category $\mathbb
%H(U,{\Gamma})$.
A commutative subalgebra ${\Gamma}\subseteq U$ is called a {\em
Harish-Chandra subalgebra\/} of $U$ \cite{dfo:hc} if for any $a\in
U$ the ${\Gamma}$-bimodule ${\Gamma} a {\Gamma}$ is finitely
generated both as left and as right ${\Gamma}$-module. The
property of $\Gamma$ to be a Harish-Chandra subalgebra is
important for the effective study of the category $ \mathbb
H(U,{\Gamma})$. In this case, for any finite-dimensional
${\Gamma}$-module $X$ the module $U\otimes_{{\Gamma}} X$ is a
Harish-Chandra module. For any $a\in U$ set \ben X_a \ = \ \{
({\mathbf m},{\mathbf n})\in \Sp {\Gamma} \times \Sp {\Gamma} \ |
\  \Gamma/{\mathbf n} \ \ \text{ is a subquotient of } \ \
{\Gamma} a {\Gamma}/ {\Gamma} a {\mathbf m} \}.
\end{equation*}
Equivalently, $({\mathbf
m},{\mathbf n})\in$ $X_a$ if and only if (${\Gamma}/ {\mathbf n})
\otimes_{{\Gamma}} {\Gamma} a {\Gamma} \otimes_{{\Gamma}}
({\Gamma}/{\mathbf m})\ne 0$. Denote by  {$\Delta$}{} the minimal
equivalence on $\Sp {\Gamma}$ containing all $X_a$, $a\in U$ and
by  {$ {\Delta}(A, {{\Gamma}})$}{} the set of the
$\Delta$-equivalence classes on $\Sp {\Gamma}$. Then for any $a\in
U$ and ${\mathbf m}\in$ $\Sp {\Gamma}$ we have
\beql{equ_1}
a M({\mathbf m}) \ \subseteq \displaystyle {\sum_{({\mathbf
m},{\mathbf n})\in X_a}}M({\mathbf n}),\qquad \mathbb
H(U,{\Gamma})= \bigoplus_{D\in \Delta(U,{\Gamma})}\mathbb
H(U,{\Gamma},D).
\end{equation}
Define a category ${\mathscr A}$ $={\mathscr A}_{U,{\Gamma}}$ with
the set of objects
$\Ob{\mathscr A}=\Sp{\Gamma}$ and with the space of morphisms
${\mathscr A}({\mathbf m},{\mathbf n})$ from
${\mathbf m}$ to ${\mathbf n}$, where
\beql{equ_2}
{\mathscr A}({\mathbf m},{\mathbf n})= {\Gamma}_{{\mathbf n}}
\otimes_{{\Gamma}} U \otimes_{{\Gamma}} {\Gamma}_{{\mathbf m}}
\end{equation}
and
$\displaystyle {\Gamma}_{{\mathbf m}}=$
$\displaystyle \lim_{\leftarrow n}{\Gamma}/{\mathbf m}^n$
is the
completion of ${\Gamma}$ by ${\mathbf m} \in$ $\Sp {\Gamma}$.
Then
the space ${\mathscr A}({\mathbf m},{\mathbf n})$ has a natural structure
of
${\Gamma}_{{\mathbf n}}{-}{\Gamma}_{{\mathbf m}}$-bimodule.
We have the decomposition
\ben
{\mathscr A}=\bigoplus_{D\in
\Delta(U,{\Gamma})}{\mathscr A}(D),
\end{equation*}
where ${\mathscr A}(D)$ is
the restriction of ${\mathscr A}$ on $D$. The category  ${\mathscr
A}$ is endowed with the topology of the inverse limit while the
category of ${\Bbbk}$-vector spaces (${\Bbbk}\dmo$) is endowed with the
discrete topology. Consider the category ${\mathscr A}\dmo_d$ of
continuous functors $M:{\mathscr A}\to{\Bbbk}\dmo$. We call them
{\em discrete modules\/}{}
following the terminology of
\cite[Section~1.5]{dfo:hc}. For any discrete
${\mathscr A}$-module $N$ define a Harish-Chandra $U$-module
\ben
\mathbb F(N)=\underset{{\mathbf m}\in\Sp {\Gamma}}
{\bigoplus}N({\mathbf m}).
\end{equation*}
Furthermore, for $x\in N({\mathbf m})$ and $a\in U$ set
\ben
ax=\sum_{{\mathbf n}\in \Sp {\Gamma}} a_{{\mathbf n}}x
\end{equation*}
where
$a_{{\mathbf n}}$ is the image of $a$ in ${\mathscr A}({\mathbf
m},{\mathbf n})$. For any morphism $f:M\to N$ in the category
${\mathscr A}\dmo_d$ set
\ben
\mathbb F(f)=\underset{{\mathbf m}\in \Sp {\Gamma}}{\bigoplus} f({{\mathbf m}}).
\end{equation*}
We thus
have a functor
 $\mathbb F:{\mathscr A}\dmo_d\to\mathbb H(U,{\Gamma})$.

\bpr\label{prop:HC-mod-equiv-A-mod} \bra\cite{dfo:hc}, Theorem 17.\ket
 The functor
$\mathbb F$ is an equivalence. \epr

 We will identify a discrete
${\mathscr A}$-module $N$ with the corresponding Harish-Chandra
module $\mathbb F(N)$.
For ${\mathbf m}\in \Sp \Gamma$ denote by $\hat{\mathbf m}$ the
completion of $\mathbf m$. Consider the two-sided ideal $I\subseteq
\mathscr A$ generated by the completions $\hat{\mathbf m}$ for all $\mathbf m\in
\Sp \Gamma$ and set ${\mathscr A}_W= {\mathscr A}/I$.
Proposition~\ref{prop:HC-mod-equiv-A-mod} implies the following.

\bco\label{cor:hc-weight-equiv-Aw}
 The categories  {$ \mathbb HW(U,{\Gamma})$}{} and
${\mathscr A}_W\dmo$ are equivalent. \eco

 The subalgebra
${\Gamma}$ is called  {\em big in\/} ${\mathbf m}\in$ $\Sp {\Gamma}$ if
${\mathscr A}({\mathbf m},{\mathbf m})$ is finitely generated as a
left (or, equivalently, right)
${\Gamma}_{{\mathbf m}}$-module.

\ble{}\label{lem:Gamma-big} \bra\cite{dfo:hc}, Corollary 19.\ket If
${\Gamma}$ is big in ${\mathbf m}\in$ $\Sp {\Gamma}$ then there
exist finitely many non-isomorphic irreducible Harish-Chandra
$U$-modules $M$ such that $M({\mathbf m})\ne 0$. For any such
module $\dim M({\mathbf m})<\infty$. \ele

 \vspace{0,2cm}

\subsection{Special PBW algebras}
Let $U$ be an associative
algebra over ${\Bbbk}$  endowed with an
increasing filtration $\{U_{i} \}_{i\in\mathbb
Z}$, $U_{-1} =$ $\{ 0 \}$, $U_{0} =$
${\Bbbk}$, $U_iU_j\subseteq U_{i+j}$.
 For $u\in U_i\setminus U_{i-1}$ set $\deg u=i$.
 Let $\overline  U=\gr \ U$  be the associated
graded algebra
\ben
{\overline  U}=
\bigoplus_{i=0}^{\infty} U_{i }/U_{i-1 }.
\end{equation*}
For $u\in U$ denote by
$\overline u$ its image in $\overline U$  and for a subset
$S\subseteq U$ set ${\overline  S}=$ $\{ {\overline s}\,|\,s\in S\}$
$\subseteq {\overline  U}$.
 The algebra $U$ is called a  \emph{special} PBW algebra if any element of
$U$ can be written uniquely as a linear combination of ordered
monomials in some fixed generators of $U$ and if  $\overline U$ is
a polynomial algebra. Such algebras were introduced  in
\cite{fo:kt}.

 Let ${\Lambda}=$
${\Bbbk}[X_1,\ldots,X_n]$ be  a polynomial algebra.
 A
sequence $g_1,\ldots$, $g_t\in{\Lambda}$ is called  {\em regular\/}
(in $\Lambda$) if the class of $g_i$ in
${\Lambda}/(g_1,\ldots,g_{i-1})$ is non-invertible and is not a
zero divisor for any $i=1,\ldots,t$.

Next proposition contains some properties of regular
sequences which will be used in the sequel; see e.g. \cite{bh:cm}.

\bpr\label{prop:reg-seq-properties}
\begin{enumerate}

 \item\label{item:reg-seq-subst-0} A
sequence of the form $X_1, \ldots, X_r, G_1, \ldots G_t$
%with $G_1, \ldots, G_t\in \Lambda$
is regular  in $\Lambda$ if and only if the
sequence $g_1, \ldots, g_t$ is regular in $\Bbbk[X_{r+1},$ $
\ldots,$ $ X_n]$, where $g_i(X_{r+1}, \ldots, X_n)=G_i(0, \ldots,
0, X_{r+1}, \ldots, X_n)$.

\item\label{item:reg-seq-fist-is-product} A sequence $g_1g_1',
g_2, \ldots, g_t$ is regular if and only if the sequences $g_1,
g_2, \ldots, g_t$ and $g_1', g_2, \ldots, g_t$ are regular.
\end{enumerate}
\epr

The following analogue of Kostant theorem \cite{k:lg} is valid
for special PBW algebras.

\bth{}\label{prop:sPBW-alg-free-over-alg-reg-seq}\cite{fo:kt}
Let $U$ be a special PBW algebra and let
 $g_1$, $\ldots$, $g_t\in U$
be mutually commuting elements such that ${\overline  g}_1$,
$\ldots$, ${\overline  g}_t$ is a regular sequence in
${\overline  U}$, and let ${\Gamma}=$ ${\Bbbk}[g_1,\ldots,g_t]$.
Then
$U$ is a free left (right) ${\Gamma}$-module. Moreover ${\Gamma}$ is a direct
 summand of $U$.
\eth

\vspace{0,5cm}

\section{Freeness of $\Y_p(\gl_2)$ over $\Gamma$}
\setcounter{equation}{0}

Let $p$ be a positive integer.
The {\it level $p$ Yangian\/} $\Y_p(\gl_2)$ for the Lie algebra
$\gl_2$ \cite{c:qg} can be defined as the algebra over $\Bbbk$
with generators $t_{ij}^{(1)}, \dots, t_{ij}^{(p)}$, $i,j=1,2$,
subject to the relations
\beql{defrelp}
[T_{ij}(u), T_{kl}(v)]=\frac{1}{u-v}(T_{kj}(u)\tss T_{il}(v)-T_{kj}(v)\tss T_{il}(u)),
\end{equation}
where $u,v$ are formal variables and
\ben
T_{ij}(u)=\delta_{ij}\ts u^p+\sum_{k=1}^p t_{ij}^{(k)}\ts u^{p-k}\in \Y_p(\gl_2)[u].
\end{equation*}
Explicitly, \eqref{defrelp} reads
\ben
[t^{(r)}_{ij}, t^{(s)}_{kl}] =\sum_{a=1}^{\min(r,s)}
\big(t^{(a-1)}_{kj} t^{(r+s-a)}_{il}-t^{(r+s-a)}_{kj} t^{(a-1)}_{il}\big),
\end{equation*}
where $t^{(0)}_{ij}=\delta_{ij}$ and $t_{ij}^{(r)}=0$ for $r\geq p+1$.
Note that the level $1$ Yangian  $\Y_1(\gl_2)$ coincides with
the universal enveloping algebra $\U(\gl_2)$.
 Set
\ben
\deg t_{ij}^{(k)}=k\qquad\text{for}\quad  i,j=1,2\quad
\text{and} \quad k=1, \ldots, p.
\end{equation*}
This defines a natural
 filtration on the Yangian $\Y_p(\gl_2)$. The corresponding graded algebra
 will be denoted by $\overline {\Y}_p(\gl_2)$.
We have the following analogue of the Poincar\'e--Birkhoff--Witt
theorem for the algebra $\Y_p(\gl_2)$.

\bpr\label{prop:pbw}\bra\cite{c:qg}; see also \cite{m:ce}\ket Given an
arbitrary linear ordering on the set of generators
$t_{ij}^{(k)}$, any element of the algebra $\Y_p(\gl_2)$ is
uniquely written as a linear combination of ordered monomials in
these generators. Moreover, the algebra $\overline {\Y}_p(\gl_2)$
is a polynomial algebra in generators
 $\overline {t}_{ij}^{(k)}$.
\epr

Proposition \ref{prop:pbw} implies that $\Y_p(\gl_2)$ is a special PBW
algebra.
Denote by $D(u)$ the {\it quantum determinant\/}
\beql{qdet}
\bal
D(u)&=T_{11}(u)\tss T_{22}(u-1)-T_{21}(u)\tss T_{12}(u-1)\\
{}&=T_{11}(u-1)\tss T_{22}(u)-T_{12}(u-1)\tss T_{21}(u)\\
{}&=T_{22}(u)\tss T_{11}(u-1)-T_{12}(u)\tss T_{21}(u-1)\\
{}&=T_{22}(u-1)\tss T_{11}(u)-T_{21}(u-1)\tss T_{12}(u).
\eal
\end{equation}
Clearly, $D(u)$ is a monic polynomial in $u$ of degree $2p$,
\beql{qdetdec}
D(u)=u^{2p}+d_1\ts u^{2p-1}+\cdots+d_{2p},\qquad d_i\in\Y_p(\gl_2).
\end{equation}
It was shown in \cite{c:ni, c:qg} (see also \cite{m:ce}
for a different proof) that the coefficients $d_1,\dots,d_{2p}$
are algebraically independent generators
of the center of the algebra $\Y_p(\gl_2)$.
Denote by $\Gamma$ the subalgebra of $\Y_p(\gl_2)$ generated by
the coefficients of $D(u)$ and by the elements
$t_{22}^{(k)}$, $k=1, \ldots, p$.
This algebra is obviously commutative. We will show later
(Corollary~\ref{thm:GT-HC}) that
 $\Gamma$ is a Harish-Chandra subalgebra in
$\Y_p(\gl_2)$.

\ble{}\label{lem:gen-of-GT-form-reg-seq} The sequence $\overline
{t}_{22}^{(1)}, \ldots, \overline {t}_{22}^{(p)}, \overline{d}_1,
\ldots, \overline{d}_{2p}$ of the images of the generators of
$\Gamma$  is regular  in $\overline {\Y}_p(\gl_2)$. \ele

\begin{proof}
Let us set
\ben
t_i=\t_{11}^{(i)}+\t_{22}^{(i)},\quad i=1,\ldots, p
\qquad\text{and}\qquad
\Delta_{i,j}=\t_{11}^{(i)}\t_{22}^{(j)}-\t_{21}^{(i)}\t_{12}^{(j)},
\quad
i,j=1, \ldots, p.
\end{equation*}
It follows from \eqref{qdetdec} that
$$\overline{D}(u)=u^{2p}+\sum_{i=1}^{2p}\overline{d}_i\ts u^{2p-i},
$$
with
\ben
\bal
\overline{d}_i&=t_i+\sum_{j=1}^{i-1}\Delta_{j,i-j}\qquad &&\text{for}\quad
i=1, \ldots, p  \qquad \text{and}\\
\overline{d}_i&=\sum_{j=i-p}^{p}\Delta_{j,i-j}\qquad &&\text{for}\quad
i=p+1,\ldots, 2p.
\eal
\end{equation*}
Hence we need to show that the sequence
$$\t_{22}^{(1)}, \ldots, \t_{22}^{(p)}, t_1, t_2+\Delta_{11}, \ldots, t_p+
\sum_{i=1}^{p-1}\Delta_{i,p-i}, \sum_{i=1}^p\Delta_{i,p+1-i},
\ldots, \Delta_{pp}$$ is regular. We will denote by $\nabla_i$ the
result of the substitution  $\t_{22}^{(1)}= \cdots=
\t_{22}^{(p)}=0$ in $\overline{d}_i$, $i=1, \ldots, 2p$. By
Proposition \ref{prop:reg-seq-properties}\ts\eqref{item:reg-seq-subst-0},
we only need to show the regularity
of the sequence
$$\nabla_{1}, \ldots,
\nabla_{2p}.$$  Consider the automorphism
$\phi$ of $\overline{\Y}_p(\gl_2)/I$ given by
\ben
\t_{11}^{(i)}\mapsto
\nabla_i,\qquad
 \t_{21}^{(i)}\mapsto \t_{21}^{(i)},\qquad \t_{12}^{(i)}\mapsto
 \t_{12}^{(i)}, \qquad \text{for}\quad
 i=1, \ldots, p,
\end{equation*}
where $I$ is the ideal generated by
$\t_{22}^{(1)}, \ldots, \t_{22}^{(p)}$. Since the regularity
of a sequence is preserved by automorphisms,
it is sufficient to demonstrate the regularity of
the sequence
 $$\t_{11}^{(1)}, \ldots, \t_{11}^{(p)}, \nabla_{p+1}, \ldots,
\nabla_{2p}.$$
Since the elements $\nabla_{i}$ do not depend on the $\t_{11}^{(k)}$,
Proposition \ref{prop:reg-seq-properties}\ts\eqref{item:reg-seq-subst-0} implies
that this is equivalent to the regularity
of the sequence $\nabla_{p+1},$ $ \ldots,$ $ \nabla_{2p}.$ For
each pair of indices $k,l\in \{1, \ldots, p\}$ and any index
$1\leq a\leq\max\{k,l\}$, consider the
sequence of $a$ elements
which occupy the rows of the table  $s(k,l,a)$ below

 $$\left(\begin{array}{l}
\phantom{\Bigl(}\t_{21}^{(k)}\t_{12}^{(l)}\\
\phantom{\Bigl(}\t_{21}^{(k-1)}\t_{12}^{(l)}+\t_{21}^{(k)}\t_{12}^{(l-1)} \\
\phantom{\Bigl(}\t_{21}^{(k-2)}\t_{12}^{(l)}+\t_{21}^{(k-1)}
\t_{12}^{(l-1)}+\t_{21}^{(k)}\t_{12}^{(l-2)}
\\
\phantom{\Bigl(}\vdots\\
\phantom{\Bigl(}\t_{21}^{(k-a+1)}\t_{12}^{(l)}+\t_{21}^{(k-a+2)}\t_{12}^{(l+1)}+
\cdots+\t_{21}^{(k)}\t_{12}^{(l-a+1)}\\
\end{array}\right)$$
Note that when $k=l=a=p$ the rows of the table are exactly the
elements $\nabla_i$, $i=p+1, \ldots, 2p$. We will show by
induction on $a$ that the sequence of rows of $s(k,l,a)$ is regular.
Note that $s(k,l,1)$ consists of the single
element $\t_{21}^{(k)}\t_{12}^{(l)}$ and
is obviously regular.  Now let $a>1$. Consider the following two
tables which we denote by $s'(k,l,a)$ and $s''(k,l,a)$, respectively.

 $$\left(\begin{array}{l}
\phantom{\Bigl(}\t_{21}^{(k)}\\
\phantom{\Bigl(}\t_{21}^{(k-1)}\t_{12}^{(l)}+\t_{21}^{(k)}\t_{12}^{(l-1)} \\
\phantom{\Bigl(}\vdots\\
\phantom{\Bigl(}\t_{21}^{(k-a+1)}\t_{12}^{(l)}+\cdots+
\t_{21}^{(k)}\t_{12}^{(l-a+1)}\\
\end{array}\right), \qquad
 \left(\begin{array}{l}
\phantom{\Bigl(}\t_{12}^{(l)}\\
\phantom{\Bigl(}\t_{21}^{(k-1)}\t_{12}^{(l)}+\t_{21}^{(k)}\t_{12}^{(l-1)} \\
\phantom{\Bigl(}\vdots\\
\phantom{\Bigl(}\t_{21}^{(k-a+1)}\t_{12}^{(l)}+\cdots+
\t_{21}^{(k)}\t_{12}^{(l-a+1)}\\
\end{array}\right)$$
Due to
Proposition~\ref{prop:reg-seq-properties}\ts\eqref{item:reg-seq-fist-is-product},
it is sufficient to verify the regularity of both $s'(k,l,a)$ and $s''(k,l,a)$.
Using again
Proposition~\ref{prop:reg-seq-properties}\ts\eqref{item:reg-seq-subst-0},
substitute $\t_{21}^{(k)}=0$ in
$s'(k,l,a)$ and $\t_{12}^{(k)}=0$ in $s''(k,l,a)$. It is easy to see
that after this substitution we obtain the tables $s(k-1, l,a-1)$ and
$s(k, l-1,a-1)$, respectively.
By the induction hypothesis, both of them
are regular and so is $s(k,l,a)$.
In particular, the sequence $s(p,p,p)$ is regular which completes
the proof.
\end{proof}

Using the regularity of the sequence $\overline
{t}_{22}^{(1)}, \ldots, \overline {t}_{22}^{(p)}, \overline{d}_1,
\ldots, \overline{d}_{2p}$  we immediately obtain the following.

\bco\label{cor:gen-center} The generators ${t}_{22}^{(1)}, \ldots,
{t}_{22}^{(p)}, {d}_1, \ldots, {d}_{2p}$ of $\Gamma$ are
algebraically independent. \eco

Combining Lemma \ref{lem:gen-of-GT-form-reg-seq} with Proposition
\ref{prop:sPBW-alg-free-over-alg-reg-seq} we obtain our first main result.

\bth\label{thm:freedom}
\begin{enumerate}
\item\label{thm:Y-free-over-Gamma} $\Y_p(\gl_2)$ is free as a
left (right) module over $\Gamma$. Moreover $\Gamma$ is a direct
summand of $\Y_p(\gl_2)$.

\item\label{thm:char-of-Gamma-extends-to-Y-mod}  Any
${\mathbf m}\in \Sp \Gamma$
 extends to an irreducible $\Y_p(\gl_2)$-module.
\end{enumerate}
\eth

For a subset $P\subseteq \Y_p(\gl_2)$ denote by $\mathbb D(P)$
the set of all $x\in \Y_p(\gl_2)$ such that there
exists $z\in{\Gamma}$, $z\ne 0$ for which $z
x\in$ $P$.

\bco\label{cor:D(P)-is-fin-gen-over-center} Let $P\subseteq
\Y_p(\gl_2)$ be a finitely generated left ${\Gamma}$-module then
$\mathbb D(P)$ is a finitely generated left ${\Gamma}$-module.
\eco

\begin{proof}
Since $\Gamma$ is a domain then $\mathbb D(P)$ is a
$\Gamma$-submodule in $\Y_p(\gl_2)$. Using the fact that
$\Y_p(\gl_2)$ is a free left ${\Gamma}$-module  we conclude that
$\Y_p(\gl_ 2)\simeq$ $F_P\oplus F$ where $F_P$ and $F$ are free
left ${\Gamma}$-modules, $F_P$ has a finite rank and $P\subseteq
F_P$. Then $\mathbb D(P)\subseteq$ $F_P$ and hence it is finitely
generated as a module over a noetherian ring.
\end{proof}

\section{Harish-Chandra modules for $\gl_2$ Yangians}\label{sec:Weight}
\setcounter{equation}{0}

In this section we introduce
universal Harish-Chandra modules $M(\ell)$.
We also describe their structure in an explicit form
in the case of generic parameters $\ell$.

Let $\L$ be a polynomial algebra in the variables $b_1, \ldots, b_p,
g_1, \ldots g_{2p}$. Define a $\Bbbk$-homomorphism $\imath:
\Gamma\to \L$ by
\beql{imath}
\imath(t_{22}^{(k)})=\sigma_{k,p}(b_1,\ldots, b_p),\qquad
\imath(d_{i})= \sigma_{i,2p}(g_1,\ldots, g_{2p}),
\end{equation}
where $\sigma_{i,j}$ is the $i$-th elementary symmetric
polynomial in $j$ variables. Due to Corollary~\ref{cor:gen-center},
$\imath$ is injective.
We will identify the elements of
$\Gamma$ with their images in $\L$ and treat them as polynomials
in the variables $b_1, \ldots, b_p, g_1, \ldots g_{2p}$ invariant
under the action of the group $S_p \times S_{2p}$ . Set $\mL=\Sp
\L$. We will identify $\mL$ with $\Bbbk^{3p}$. If
\ben
\beta=(\beta_1,\ldots, \beta_p),\qquad\gamma=(\gamma_1,\ldots,
\gamma_{2p})\qquad \text{and}\qquad
\ell=(\beta_1,\ldots, \beta_p, \gamma_1,\ldots,
\gamma_{2p})
\end{equation*}
then we shall write $\ell=(\beta, \gamma)$.  The monomorphism
$\imath$ induces the epimorphism
\beql{imathstar}
\imath^*: \mL\rightarrow \Sp
\Gamma.
\end{equation}
If $\ell\in \mL$ and $\bm=\imath^*(\ell)$ then $D(\ell)$
will denote  the equivalence class of $\bm$ in
$\Delta(\Y_p(\gl_2),\Gamma)$; see Section~\ref{subsec:Harish}.

Let $\mP_0\subseteq \mL$, $\mP_0\simeq \Z^{p}$, be the lattice
generated by the elements $\delta_i\in \Bbbk^{3p}$ for  $i=1, \ldots, p$, where
\ben
\delta_i=(\delta^1_i, \ldots, \delta_i^{3p}),
\qquad\delta^j_i=\delta_{ij}.
\end{equation*}
Then $\mP_0$   acts
on $\mL$ by shifts $\delta_i(\ell):=\ell+\delta_i$. Furthermore, the
group $S_p \times S_{2p}$ acts on $\mL$ by permutations. Thus  the
semidirect product $\bW$  of the groups $S_p \times S_{2p}$ and
$\mP_0$ acts on $\mL$ and $\L$. Denote by $S$ a multiplicative set
in $\L$ generated by the elements $b_i-b_j-m$ for all $i\neq j$
and all $m\in \Z$ and by $\bL$ the localization of $\L$ by $S$.
Note that $S$ is invariant under the action of $\bW$ and hence
$\bW$ acts on $\bL$ as well.

For arbitrary $3p$-tuple $\ell=(\beta, \gamma)\in \mL$ set
\ben
\beta(u)=(u+\beta_1)\cdots (u+\beta_p),\qquad
\gamma(u)=(u+\gamma_1)\cdots (u+\gamma_{2p}).
\end{equation*}
Let $I_{\ell}$ be the left
ideal of $\Y_p(\gl_2)$ generated by the coefficients of the
polynomials $T_{22}(u)-\be(u)$ and $D(u)-\ga(u)$. Define the
corresponding quotient module over $\Y_p(\gl_2)$ by
\beql{mlmodule}
M(\ell)=\Y_p(\gl_2)/I_{\ell}.
\end{equation}
We shall call it the {\it universal module\/}.
It follows from Theorem~\ref{thm:freedom} that $I_{\ell}$ is
a proper ideal of $\Y_p(\gl_2)$ and so $M(\ell)$ is a non-trivial
module. It is clear that if
$V$ is an arbitrary Harish-Chandra $\Y_p(\gl_2)$-module generated by
a non-zero $\eta\in V$ such that $D(u)\eta=\gamma(u)\eta$ and
$T_{22}(u)\eta=\beta(u)\eta$
 then $V$ is a
homomorphic image of $M(\ell)$.

Set $\mP_1=$ $\Sp \bL$ $\subseteq {\mL}$, i.e. $\mP_1$ consists of
\emph{generic} $3p$-tuples $\ell=(\beta,\gamma)$
 such that
\beql{genbeta}
\beta_i-\beta_j\not\in \Z\qquad \text{for all}\quad i\neq j.
\end{equation}
 If $\ell\in\mP_1$ then the modules
 from the category $\mathbb H(\Y_p(\gl_2),$ ${\Gamma},$ $D(\ell))$
 are called
\emph{generic} Harish-Chandra modules.

\subsection{Weight modules}
For $\ell=(\beta, \gamma)\in \mL$
 the category $\mathbb
HW(\Y_p(\gl_2),$ ${\Gamma},$ $D(\ell))$ consists of finitely
generated weight modules $V$ with central character $\ga$ and with
$\Supp V\subseteq  D(\ell)$. We shall denote this category by
$R_\ell$ for brevity. If $\ell\in \mP_1$ then the modules from $R_\ell$  will
be called \emph{generic} weight modules.

A $\Y_p(\gl_2)$-module $V$ is an object of $R_{\ell}$  if $V$ is a
direct sum of its \emph{weight} subspaces:
\ben
V=\underset{\ell\in \mL}{\bigoplus}\ts V_{\ell}, \qquad
V_{\ell}=\{\eta\in V\ |\ T_{22}(u)\tss\eta=\be(u)\tss\eta,\quad
D(u)\tss\eta=\ga(u)\tss\eta\}.
\end{equation*}

If $V\in  R_{\ell}$ then
 we shall simply write $V_{\be}$ instead of
$V_{\ell}$ and identify $\Supp V$ with the set of all $\be$ such
that the subspace $V_{\be}$ is nonzero. The next lemma describes
the action of the Yangian generators on the weight subspaces;
cf. \eqref{equ_1}.

\ble\label{lem:wgen}
Let $V$ be a
generic weight $\Y_p(\gl_2)$-module and let $\be=(\beta_1, \ldots,
\beta_p)\in \Supp V$. Then
\beql{t2112}
T_{21}(u)V_{\be}\subseteq
\sum_{i=1}^p V_{\be+\delta_i}\qquad\text{and} \qquad
T_{12}(u)V_{\be}\subseteq \sum_{i=1}^p V_{\be-\delta_i}
\end{equation}
where $\be\pm \delta_i = (\beta_1, \ldots, \beta_i\pm 1, \ldots ,
\beta_p)$. \ele

\begin{proof}
First we show that $T_{21}(-\beta_i)V_{\be}\subseteq V_{\be+\delta_i}$
for all $i=1, \ldots, p$. Since
$$
T_{22}(u-1)T_{21}(u)=T_{21}(u-1)T_{22}(u)
$$
we have
$$
T_{22}(-\beta_i-1)T_{21}(-\beta_i)\ts\eta=T_{21}(-\beta_i-1)T_{22}(-\beta_i)\ts\eta=0
$$
for all $\ts\eta\in V_{\be}$. Also,
\ben
\bal
T_{22}(-\beta_j)T_{21}(-\beta_i)\ts\eta&=
(\beta_i-\beta_j)^{-1}(T_{21}(-\beta_i)
T_{22}(-\beta_j)-T_{21}(-\beta_j)T_{22}(-\beta_i))\ts\eta\\
{}&+T_{21}(-\beta_i)T_{22}(-\beta_j)\ts\eta=0
\eal
\end{equation*}
since $T_{22}(-\beta_k)\ts\eta=0$ for all $k=1, \ldots, p$.
Using the fact that $\beta_i-\beta_j \notin \Z$ we conclude that
$T_{21}(-\beta_i)V_{\be}\subseteq V_{\be+\delta_i}$
for all $i=1, \ldots, p$. Since $T_{21}(u)$
is a polynomial of degree $p-1$ in $u$
and $\beta_i\neq \beta_j$ if $i\neq j$, we thus get the first containment
of \eqref{t2112}. The second is verified in the same way with the use
of the identity
$T_{22}(u)T_{12}(u-1)=T_{12}(u)T_{22}(u-1)$.
\end{proof}

\bco\label{cor:suppind} If $V$ is indecomposable generic weight
module over $\Y_p(\gl_2)$ and $\be\in \Supp V$ then $\Supp
V\subseteq \be +\Z^p$. \qed\eco

\ble\label{lem:comr} If $V$ is a generic weight
$\Y_p(\gl_2)$-module with the central character $\ga$ then for any
$\be=(\beta_1, \ldots, \beta_p)\in \Supp V$ and any $\ts\eta\in
V_{\be}$ we have
$$
T_{12}(-\beta_r)T_{21}(-\beta_s)\ts\eta=T_{21}(-\beta_s)T_{12}(-\beta_r)\ts\eta,
$$
if $s\neq r$, and
\ben
\bal
T_{12}(-\beta_i-1)T_{21}(-\beta_i)\ts\eta&=-\ga(-\beta_i)\ts\eta,\\
T_{21}(-\beta_i+1)T_{12}(-\beta_i)\ts\eta&=-\ga(-\beta_i+1)\ts\eta.
\eal
\end{equation*}
\ele

\begin{proof}
The first equality follows from the defining relations
\eqref{defrel}. The two remaining follow from \eqref{qdet}.
\end{proof}

The following corollary is immediate from Lemma~\ref{lem:comr}.

\bco\label{cor:gasupp} Let $V$ be a generic weight
$\Y_p(\gl_2)$-module with the central character $\ga$ and let
$\be=(\beta_1, \ldots, \beta_p)\in \Supp V$.
\begin{enumerate}
\item If $\ga(-\beta_i)\neq 0$ then
$\Ker T_{21}(-\beta_i)\cap V_{\be}=0$.
\item  If $\ga(-\beta_i +1)\neq 0$
then $\Ker T_{12}(-\beta_i)\cap V_{\be}=0$.
\item  If $V$ is indecomposable
and $\ga(-\beta_i +k)\neq 0$ for all $k\in \Z$ then
$$
\Ker T_{21}(-\psi_i)\cap V_{\psi}=
\Ker T_{12}(-\psi_i)\cap V_{\psi}=0
$$
for all $\psi=(\psi_1, \ldots, \psi_p)\in \Supp V$.     \qed
\end{enumerate}
\eco

Since the universal module $M(\ell)$ is non-trivial,
the image of $1$ in $M(\ell)$ is nonzero. We
shall denote this image by $\xi$.
Assume that $\be$ satisfies
the genericity condition \eqref{genbeta}.
For any $(k)=(k_1,\dots,k_p)\in\Z^p$ define the corresponding vector
of the module $M(\ell)$ by
\beql{xikold}
\bal
\xi^{(k)}=&\prod_{i,\ k_i> 0} T_{21}(-\be_i-k_i+1)\cdots T_{21}(-\be_i-1)\ts
T_{21}(-\be_i)\\
{}\times &\prod_{i,\ k_i< 0} T_{12}(-\be_i-k_i-1)\cdots T_{12}(-\be_i+1)\ts
T_{12}(-\be_i)\ts \xi.
\eal
\end{equation}

\bth\label{thm:basis} The vectors $\xi^{(k)}, (k)\in \Z^p$ form a
basis of $M(\ell)$. Moreover, we have the formulas
\beql{T22act}
T_{22}(u)\ts\xi^{(k)}=\prod_{i=1}^p(u+\beta_i+k_i)\ts\xi^{(k)},
\end{equation}
\beql{T21u12u}
\bal
T_{21}(u)\ts \xi^{(k)}&=\sum_{i=1}^p A_i(k)\ts
\frac{(u+\be_1+k_1)\cdots\wedge_i\cdots (u+\be_p+k_p)}
{(\be_1-\be_i+k_1-k_i) \cdots\wedge_i\cdots (\be_p-\be_i+k_p-k_i)}
\ts \xi^{(k+\delta_i)}, \\
T_{12}(u)\ts \xi^{(k)}&=\sum_{i=1}^p B_i(k) \ts
\frac{(u+\be_1+k_1)\cdots\wedge_i\cdots (u+\be_p+k_p)}
{(\be_1-\be_i+k_1-k_i) \cdots\wedge_i\cdots (\be_p-\be_i+k_p-k_i)}
\ts \xi^{(k-\delta_i)},
\eal
\end{equation}
where
$$
A_i(k)=
\begin{cases} 1\quad &\text{if}\quad     k_i\geq 0\\
-\ga(-\be_i-k_i)\quad &\text{if}\quad    k_i< 0
\end{cases}
$$
and
$$
B_i(k)=
\begin{cases} -\ga(-\be_i-k_i+1)\quad &\text{if}\quad    k_i> 0\\
1\quad &\text{if}\quad   k_i\leq 0.
\end{cases}
$$
The action of $T_{11}(u)$ is found from the relation
\beql{t11u}
\Big(T_{11}(u)\tss T_{22}(u-1)-T_{21}(u)\tss T_{12}(u-1) \Big)\ts \xi^{(k)}
=\ga(u)\ts \xi^{(k)}.
\end{equation}
\eth

\begin{proof}
We start by proving the formulas for the action of the generators
of $\Y_p(\gl_2)$.
Formula \eqref{T22act} follows by induction
with the use of the relations
\begin{align}\label{T22T21}
T_{22}(u)\ts T_{21}(v)&=\frac{u-v+1}{u-v}\ts  T_{21}(v)\ts T_{22}(u)-
\frac{1}{u-v}\ts  T_{21}(u)\ts T_{22}(v)\\
\intertext{and}
\label{T22T12}
T_{22}(u)\ts T_{12}(v)&=\frac{u-v-1}{u-v}\ts T_{12}(v)\ts T_{22}(u)+
\frac{1}{u-v}\ts  T_{12}(u)\ts T_{22}(v)
\end{align}
implied by
\eqref{defrelp}.
By Lemma~\ref{lem:comr} we have:
if $k_i> 0$ then
\beql{k>0}
\bal
T_{21}(-\be_i-k_i)\ts\xi^{(k)}&=\xi^{(k+\delta_i)}, \\
T_{12}(-\be_i-k_i)\ts\xi^{(k)}&=-\ga(-\be_i-k_i+1)\ts\xi^{(k-\delta_i)};
\eal
\end{equation}
if $k_i< 0$ then
\beql{k<0}
\bal
T_{12}(-\be_i-k_i)\ts\xi^{(k)}&=\xi^{(k-\delta_i)}, \\
T_{21}(-\be_i-k_i)\ts\xi^{(k)}&=-\ga(-\be_i-k_i)\ts\xi^{(k+\delta_i)};
\eal
\end{equation}
and if $k_i= 0$ then
\beql{k=0}
\bal
T_{12}(-\be_i)\ts\xi^{(k)}&=\xi^{(k-\delta_i)}, \\
T_{21}(-\be_i)\ts\xi^{(k)}&=\xi^{(k+\delta_i)}.
\eal
\end{equation}
Applying the Lagrange interpolation formula we obtain
the remaining formulas.

It is implied by the formulas, that the module $M(\ell)$
is spanned by the vectors $\xi^{(k)}$.
By
\eqref{T22act} and the genericity assumption, the
$\xi^{(k)}$ are eigenvectors for $T_{22}(u)$ with
distinct eigenvalues. In order to verify their
linear independence, suppose first that $\ga(u)$
satisfies the condition
\beql{genegacond}
\ga(-\be_i-k)\ne 0\qquad\text{for all}\quad k\in\Z
\quad\text{and all}\quad i.
\end{equation}
In this case the linear independence of the $\xi^{(k)}$ follows
from the fact that each of them is nonzero.
This is implied by \eqref{k>0}--\eqref{k<0}
since $\xi\ne 0$ in $M(\ell)$.

In the case of general $\ga(u)$ let us define a $\Y_p(\gl_2)$-module
$\wt M(\ell)$ as follows. As a vector space,
$\wt M(\ell)$ is the $\Bbbk$-linear span of the basis vectors
$\wt\xi^{(k)}$ with $(k)$ running over $\Z^p$
and the action of $\Y_p(\gl_2)$ is given by
the formulas \eqref{T22act}--\eqref{t11u}, where
the $\xi^{(k)}$ should be replaced with $\wt\xi^{(k)}$.
We have to verify that the operators $T_{ij}(u)$ do satisfy
the Yangian defining relations \eqref{defrelp}.
However, the application of both sides of
\eqref{defrelp} to a basis vector $\wt\xi^{(k)}$
amounts to polynomial relations on the coefficients of $\ga(u)$.
By the previous argument, if $\ga(u)$ satisfies \eqref{genegacond}
then these relations are identities. Therefore, these identities remain
hold for an arbitrary $\ga(u)$ and thus $\wt M(\ell)$ is well defined.

Finally, consider the $\Y_p(\gl_2)$-module homomorphism
\ben
\varphi: \Y_p(\gl_2)\to \wt M(\ell),\qquad 1\mapsto \wt\xi^{(0)}.
\end{equation*}
Obviously, the ideal $I_{\ell}$ is contained in the kernel $\Ker\varphi$
and so, this defines a homomorphism  $M(\ell)\to\wt M(\ell)$
which takes $\xi^{(k)}$ to the corresponding
vector $\wt\xi^{(k)}$. Since the vectors $\wt\xi^{(k)}$
form a basis of $\wt M(\ell)$ this proves that the vectors
$\xi^{(k)}$ are linearly independent.
\end{proof}

Let us fix a $p$-tuple $\be$    satisfying the genericity condition
\eqref{genbeta} and introduce the elements
of $\Y_p(\gl_2)$
by
\ben
\bal
\tau^{(k)}=&\prod_{i,\ k_i> 0} T_{21}(-\be_i-k_i+1)\cdots T_{21}(-\be_i-1)\ts
T_{21}(-\be_i)\\
{}\times &\prod_{i,\ k_i< 0} T_{12}(-\be_i-k_i-1)\cdots T_{12}(-\be_i+1)\ts
T_{12}(-\be_i),
\eal
\end{equation*}
where $(k)$ runs over $\Z^p$.

\bco\label{cor:linindep}
The elements $\tau^{(k)}$ are linearly independent
over $\Ga$ in the right $\Ga$-module $\Y_p(\gl_2)$.
\eco

\begin{proof}
Suppose that a linear combination of the elements $\tau^{(k)}$
with coefficients in $\Ga$ is zero:
\beql{combzero}
\sum_{(k)}\tau^{(k)}\ts c_{(k)}=0,\qquad c_{(k)}\in\Ga.
\end{equation}
Apply the left hand side to the vector $\xi$ in a module $M(\ell)$
with $\ell=(\beta,\ga)$ satisfying the assumptions of Theorem~\ref{thm:basis}.
We get the relation
\ben
\sum_{(k)}c_{(k)}(\ell)\ts\xi^{(k)} =0,
\end{equation*}
where $c_{(k)}(\ell)$ is the evaluation of the polynomial $c_{(k)}$
at $T_{22}(u)=\beta(u)$ and $D(u)=\gamma(u)$.
Since the vectors $\xi^{(k)}$ form a basis of $M(\ell)$ this implies that
$c_{(k)}(\ell)=0$ for any choice of the parameters $\ga$.
Therefore, each $c_{(k)}$ does not depend on the generators $d_i$
and so it is a polynomial in the $t_{22}^{(i)}$.
However, due to the Poincar\'e--Birkhoff--Witt
theorem for the algebra $\Y_p(\gl_2)$ (Proposition~\ref{prop:pbw}),
a nontrivial relation \eqref{combzero} can only hold
if the elements $\tau^{(k)}$ are linearly dependent over $\Bbbk$.
But this is not the case because the vectors $\xi^{(k)}=\tau^{(k)}\ts\xi$
are linearly independent in $M(\ell)$ by Theorem~\ref{thm:basis}.
\end{proof}

\bre\label{rem:leftmod}
One can also produce a family
of $\Ga$-linearly independent elements for
the left $\Ga$-module $\Y_p(\gl_2)$.
They can be obtained as images of the $\tau^{(k)}$ under the
anti-automorphism
of the algebra $\Y_p(\gl_2)$ given by
\beql{antiauto}
t_{ij}^{(r)}\mapsto t_{ji}^{(r)}.
\end{equation}
For the proof we observe that
every generator of $\Ga$ is stable under this anti-automorphism.
With the exception of the case $p=1$, the elements $\tau^{(k)}$ do not
apparently constitute a basis of
$\Y_p(\gl_2)$ as a right $\Ga$-module.
\ere

\bre Given two monic polynomials $\al(u)$ and $\be(u)$ of degree
$p$ define the corresponding {\it Verma module\/}
$V(\al(u),\be(u))$ as the quotient of $\Y_p(\gl_2)$ by the left
ideal generated by the coefficients of the polynomials
$T_{11}(u)-\al(u)$, $T_{22}(u)-\be(u)$ and $T_{12}(u)$; cf.
\cite{t:sq, t:im}. Then the same argument as above shows that
$V(\al(u),\be(u))$  has a basis $\{\xi^{(k)}\}$ parameterized by
$p$-tuples of nonnegative integers $(k)$.
The formulas of Theorem~\ref{thm:basis} hold for the basis vectors
$\xi^{(k)}$, where $\ga(u)$ should be taken to be
$\al(u)\ts\be(u-1)$ which defines the central character $\gamma$
of $V(\al(u),\be(u))$. In fact, $V(\al(u),\be(u))$ is isomorphic
to the quotient of the corresponding universal module $M(\ell)$,
$\ell=(\be,\ga)$ by the submodule spanned by the vectors
$\{\xi^{(k)}\}$ such that $(k)$ contains at least one negative
component $k_i$. \ere

\bco\label{cor:suppuni} Let $\ell=(\be,\ga)\in \mP_1$.
\begin{enumerate}
\item\label{item:suppuni-onedim} The module $M(\ell)$ is a generic
weight $\Y_p(\gl_2)$-module with central character $\ga$, $\Supp
M(\ell)=\Z^p$ and all weight spaces are $1$-dimensional.
\item\label{item:suppuni-unique-max} The  module $M(\ell)$ has a
unique maximal submodule and hence a unique irreducible quotient.
\item\label{item:suppuni-and-D} The equivalence class $D(\ell)$
coincides with the set $\ell+\mP_0$.
\end{enumerate}
\eco

\begin{proof}
Statement \eqref{item:suppuni-onedim} follows immediately from
Theorem~\ref{thm:basis}. The sum of
all proper submodules of $M(\ell)$  is again a proper submodule implying
\eqref{item:suppuni-unique-max}. The statement
\eqref{item:suppuni-and-D} follows immediately from
\eqref{item:suppuni-onedim}.
\end{proof}

We will denote the unique irreducible quotient of $M(\ell)$ by
$L(\ell)$. It follows from Corollary~\ref{cor:suppuni} that all
weight spaces of $L(\ell)$ are $1$-dimensional. We can now describe all irreducible
generic weight $\Y_p(\gl_2)$-modules.

\bco\label{cor:gen} Let $\ell=(\beta,\gamma)\in \mP_1$.
\begin{enumerate}
\item There exists an irreducible generic
weight $\Y_p(\gl_2)$-module $L(\ell)$ with $L(\ell)_{\be}\neq 0$
and with central character $\ga$. Moreover, $\dim
L(\ell)_{\psi}=1$ for all $\psi \in \Supp L(\ell)$.
\item Any irreducible weight module over
$\Y_p(\gl_2)$ with central character $\ga$ generated by a nonzero
vector of  weight $\be$ is isomorphic to $L(\ell)$.
\end{enumerate}
\eco

\section{$\Gamma$ is a Harish-Chandra subalgebra}\label{sec:Some}
\setcounter{equation}{0}

In this section we adapt the results from \cite{dfo:hc} and
\cite{o:fs} for
 the Yangians. In particular, we  show that
$\Gamma$ is a Harish-Chandra subalgebra.

We have
 the following  analogue of the Harish-Chandra theorem
 for Lie algebras \cite{d:ae}.

 \bpr\label{cor:annih-all-generic=0}
 Let $x\in \Y_p(\gl_2)$ be such that
 $xM(\ell)=0$ for any
${\ell}\in
 \mP_1$.  Then $x=0$.
\epr

\begin{proof}
 Since
$M(\ell)=\Y_p(\gl_2)/I_{\ell}$  it will be sufficient to show that
the intersection $\bigcap_{\ell} I_{\ell}$ over all
$\ell\in
\mP_1$ is zero. By Theorem \ref{thm:freedom}\ts\eqref{thm:Y-free-over-Gamma},
the Yangian $\Y_p(\gl_2)$ is free as
a right module over $\Gamma$. Let $x_i, i\in \mathcal I$ be a
basis of $\Y_p(\gl_2)$ over $\Gamma$. If $x=\sum_{i\in \mathcal
I}x_iz_i$ for some $z_i\in \Gamma$ then $x\in I_{\ell}$ if and
only if $z_i(\ell)=0$ for all $i\in \mathcal I$. Since $\mP_1$ is
dense in $\mL$ in Zariski topology it follows immediately that if
$x\in \bigcap_{\ell}I_{\ell}$ then $z_i=0$ for all $i\in
\mathcal I$ and thus $x=0$. This completes the proof.
\end{proof}

 For any $\ell_0\in \mP_1$ the module $M(\ell_0)$ has
a basis $\xi^{(k)}$, $(k)\in \Z^p$ with the action of generators
of $\Y(\gl_2)$ defined
 by formulas \eqref{T22act}--\eqref{t11u}.  We will relabel the basis elements of
$M(\ell_0)$ as $\xi_{\ell}$, $\ell\in \ell_0+\mP_0$. It follows
from Theorem~\ref{thm:basis} that for every $x\in \Y_p(\gl_2)$
there exists a finite subset  $\mathcal{L}_x\subseteq \mP_0$
consisting of elements $\delta$ such that

\beql{xiell}
  x\ts \xi_{\ell}=\sum_{\delta\in \mathcal{L}_x}
\theta(x,\ell,\delta)\ts\xi_{\ell+\delta}
\end{equation}
for some nonzero coefficients $\theta(x,\ell,\delta)\in\Bbbk$. We can also
regard these coefficients as the values of the rational functions
${\theta}(x, \mathsf b, \delta)\in\bL$
at $\mathsf b=\ell$, where
$\mathsf b=(b_1, \ldots, b_p,g_1,\ldots,g_{2p})$.
Clearly, the set $\mathcal{L}_x$ is $S_p\times S_{2p}$-invariant.
%menyayutsya k_i, a znachit delta_i ,a beta_i fiksirovani
Note that for a given $x$ this formula does not depend on $\ell_0$.

We identify the
$({\Gamma}{-}{\Gamma})$-bimodule structure on $\Y_p(\gl_2)$ with the
corresponding ${\Gamma}\otimes{\Gamma}$-module
structure.
 For any $z\in \Gamma$ and any finite $S \subseteq \mL$  introduce the following
polynomial

$$F_{S,z}=F_{S,z}(z,\mathsf b)=\displaystyle
\prod_{\delta\in S} (z\otimes 1-1\otimes z(\mathsf b
+\delta))=\sum_{i=0}^{|S|}z^{i}\otimes a_i,\qquad a_i\in \L.$$

\bpr\label{prop:prop-sootn}\bra cf. \cite[Lemma 25]{dfo:hc}.\ket Let $S$ be a
finite $S_p\times S_{2p}$-invariant subset in $\mL$, $q=|S|$,  $z\in \Gamma$
and
$F_{S,z}=\displaystyle\sum_{i=0}^{q}z^{i}\otimes a_i, \quad a_i\in
\L.$ Then
\begin{enumerate}
\item\label{item:coeff-are-in-Gamma} $a_i\in \Gamma$,
$i=0,\dots, q$.

\item\label{item:polynom-annihilate}  For any $x\in \Y_p(\gl_2)$
such that $\mL_x\subseteq S$ we have
$\displaystyle\sum_{i=0}^qz^{i}xa_i=0.$
\end{enumerate}
\epr

\begin{proof} Since $S$ is $S_p\times S_{2p}$-invariant, the coefficients of the
polynomial $F_{S,z}$ are $S_p\times S_{2p}$-invariant and hence
belong to $\Gamma$ which proves \eqref{item:coeff-are-in-Gamma}.
It is sufficient to check the statement
(\ref{item:polynom-annihilate}) for $S=\mL_x$ since
$F_{S,z}=F_{S\setminus \mL_x,z}F_{\mL_x,z}$.
Let $\ell\in \mP_1$ and let  $\xi_{\ell}$ be  a basis element of
$M(\ell)$. Then
$$\sum_{i=0}^qz^{i}xa_i(\xi_{\ell})=\sum_{i=0}^qz^{i}xa_i(\ell)(\xi_{\ell})$$
$$=\sum_{i=0}^qz^{i}a_i(\ell)\sum_{\delta\in \mL_{x}}\theta(x, \ell,
\delta)\xi_{\ell+\delta}$$
$$=\sum_{\delta\in \mL_{x}}\theta(x, \ell,
\delta)\sum_{i=0}^qa_i(\ell)(z^{i}\xi_{\ell+\delta})$$
$$=\sum_{\delta\in \mL_{x}}\theta(x, \ell,
\delta)\sum_{i=0}^qa_i(\ell)z(\ell+\delta)^{i}\xi_{\ell+\delta}=\sum_{\delta\in
\mL_x}\theta(x, \ell, \delta)F_{\mL_x,z}(z(\ell+\delta),
\ell)\xi_{\ell+\delta}=0$$ since $F_{\mL_ x,z}(z(\ell+\delta),
\ell)=0$ for every $\delta \in \mL_x$. Applying Corollary
\ref{cor:annih-all-generic=0} we obtain the statement of the
proposition.
\end{proof}

The main result of this section is the following theorem.

\bth\label{thm:GT-HC}
 $\Gamma$ is a Harish-Chandra subalgebra of $\Y_p(\gl_2)$.
\eth

\begin{proof}
Following \cite[Proposition~8]{dfo:hc}, it is sufficient to show
that a $\Gamma$-bimodule  $\Gamma\, t_{ij}^{(k)}\,\Gamma$ is
finitely generated both as left and as right module for every
possible choice of indices $i,j,k$. This is obvious for $i=j=2$
since $t_{22}^{(k)}\in \Gamma$. Suppose now that  $i=2, j=1$.
We have
$d_i\ts t_{21}^{(k)}=t_{21}^{(k)}d_i$ by the centrality of $d_i$.
It follows from \eqref{T21u12u}
that $\mathcal{L }_{t_{21}^{(k)}}=\{\delta_{i}\ |\
i=1, \ldots, p\}$.  Then for $x=t_{21}^{(k)}$ we have
$$
F_{\mL_{x},t_{22}^{(i)}}=z^p\otimes
1+\sum_{l=0}^{p-1}z^{l}\otimes a_l,\qquad a_l\in \Gamma
$$
and
\beql{eq-GZ}
({t_{22}^{(i)}})^p
t_{21}^{(k)}+\sum_{l=0}^{p-1}({t_{22}^{(i)}})^{l}t_{21}^{(k)}a_l=0
\end{equation}
by Proposition \ref{prop:prop-sootn}\ts\eqref{item:polynom-annihilate}.
Hence the elements
\ben
\prod_{i=1}^{p}(t_{22}^{(i)})^{k_i}\ts t_{21}^{(k)}, \qquad 0\leq k_i<p
\end{equation*}
are generators of $\Gamma t_{21}^{(k)}\Gamma$ as a right
$\Gamma$-module.
% To show that  $\Gamma t_{21}^{(k)}\Gamma$ is finitely
%generated as a left $\Gamma$-module we consider
%$$a_0=\prod_{l=1}^p t_{22}^{(i)}(\mathsf b +\delta_l)$$ in
%\ref{eq-GZ}. Due to the theory of symmetric polynomials
%$a_0=(t_{22}^{(i)})^p+\sum_{l=0}^{p-1}(t_{22}^{(i)})^{l}y_l$ where
%$y_l$ is a polynomial in $t_{22}^{(j)}$ with $j<i$. Also note that
%$a_s$, $0<s<p$, at most contains the power $(t_{22}^{(i)})^{p-s}$.
% Hence, as above the elements
%$t_{21}^{(k)}(\prod_{i=1}^{p}(t_{22}^{(i)})^{k_i})$, $0\leq k_i<p$
%form the generators of $\Gamma\, t_{21}^{(k)}\,\Gamma$ as a left
%$\Gamma$-module.
The cases $i=1,j=2$ and $i=j=1$ are treated similarly since
\ben
\mathcal{L}_{t_{12}^{(k)}}=\{-\delta_{i}\ |\  i=1, \ldots, p\}\qquad \text{and}
\qquad
\mathcal{L}_{t_{11}^{(k)}}=\{\ts\delta_{i}-\delta_{j}\ |\  i,j=1,
\ldots, p\}.
\end{equation*}
Thus, $\Gamma t_{ij}^{(k)}\Gamma$ is finitely
generated as a right $\Gamma$-module.
The claim for the left module is proved by the application of
the anti-automorphism of the algebra $\Y_p(\gl_2)$
defined in \eqref{antiauto} where we note that every element of
$\Ga$ is stable under this anti-automorphism.
\end{proof}

\bex\label{example:relp2}
We give an explicit form of the relation \eqref{eq-GZ} for the
particular case
$i=k=p=2$. It reads
\ben
{t_{22}^{(2)}}^2\ts t_{21}^{(2)}-t_{22}^{(2)}\ts t_{21}^{(2)}
\Big(2\ts t_{22}^{(2)}+t_{22}^{(1)}\Big)+
t_{21}^{(2)}\Big( {t_{22}^{(2)}}^2+
t_{22}^{(2)}\ts t_{22}^{(1)}+t_{22}^{(2)}\Big)=0.
\end{equation*}
\eex

\section{Universal representation of the Yangian}

We will denote by $K(\Gamma)$  the field of fractions of $\Gamma$.

Let $\M_{\mP_0}(\bL)$ be the ring of locally finite matrices over $\bL$
(with a finite number of non-zero elements in each row and each column)
 with the  entries indexed by the elements of
$\mP_0$. Any $\ell\in$ $\mP_1$ defines the evaluation
homomorphism $\chi_\ell:\mathbb L\to\Bbbk$,
which induces the homomorphism of matrix algebras
$\Phi(\ell):$ ${\mathrm M}_{\mP_0}(\mathbb
L)$ ${\to}$ ${\mathrm M}_{\mP_0}({\Bbbk})$.
For $\ell, \ell'\in \mP_0$ denote by $e_{\ell \, \ell'}$ the
corresponding matrix unit in $\M_{\mP_0}(\bL)$. The group $\bW$
acts  on $\M_{\mP_0}(\bL)$ by the rule: if
$
X=(X_{\ell \,\ell'})_{\ell , \ell'\in \mP_0}
$
then
\beql{waction}
(w^{-1}\cdot X)_{\ell,\ell'}=w^{-1}\cdot
X_{w(\ell)w(\ell')}\qquad \text{for}\quad w\in \bW.
\end{equation}

Define the map
$$\G:\Y_p(\gl_2)\to \M_{\mP_0}(\bL)$$
such that for any $x\in \Y_p(\gl_2)$ and any $\ell\in \mP_0$,
$\G(x)_{\ell \, \ell'}={\theta}(x, \mathsf b + \ell, \delta)$ if
$\ell' - \ell=\delta$
 and $0$ otherwise; see \eqref{xiell}.

\ble\label{lem:W-invariants}
\begin{enumerate}
\item\label{item:G-is-repr-of-Y} $\G$ is a representation of
$\Y_p(\gl_2)$.

\item\label{item:Im-G-is-W-invar} $\G(x)$ is $\bW$-invariant for
any $x\in \Y_p(\gl_2)$. In particular,
 $\G(x)_{\o\, \o}\in K(\Gamma)$.

\item\label{item:diag-of-G(Gamma)} If $x=x(b_1, \ldots, b_p, g_1,
\ldots, g_{2p})\in \Gamma$ and $\ell=(l_1,\ldots,
l_p, 0, \ldots, 0)\in \mP_0$
then
\ben
\G(x)_{\ell \ell}=x(b_1+l_1,
\ldots, b_p+l_p, g_1, \ldots, g_{2p}).
\end{equation*}

\item\label{item:diag-of-G(Gamma)-W-inw-defined-by-0} $\G(\Gamma)$
consists of $\bW$-invariant diagonal matrices $X$ such that
$X_{\o\, \o}\in \Gamma$. In particular,
any such matrix $X$ is determined by
$X_{\o\, \o}\in \Gamma$.
\end{enumerate}
\ele

\begin{proof}
Let $\T$ be the free associative algebra with generators
$\wt t_{ij}^{\ts(k)}$, where $i,j=1,2$ and $k=1, \ldots, p$,
and let
\ben
\pi:\T\to
\Y_p(\gl_2),\qquad \wt t_{ij}^{\ts(k)}\mapsto t_{ij}^{(k)},
\end{equation*}
be the
canonical projection. Define a homomorphism $g:\T\to
\M_{\mP_0}(\bL)$ by $g(\wt t_{ij}^{\ts(k)})=\G(t_{ij}^{(k)})$ for all
suitable $i,j,k$.  To prove \eqref{item:G-is-repr-of-Y} it is
sufficient to show
 that $g(\Ker \pi)=0$. Suppose that $f\in \Ker
\pi$. Then
 $\Phi(\ell)(g(f))=0$  and thus $g(f)_{\ell'
\ell''}(\ell)=0$ for any $\ell\in \mP_1$. Since $\mP_1$ is dense
in $\Sp L$ we conclude that $g(f)=0$ implying
\eqref{item:G-is-repr-of-Y}.
The image of $\G$ is $\bW$-invariant since this holds for the
generators of $\Y_p(\gl_2)$; see \eqref{T22act}--\eqref{t11u}. For any
$\sigma\in S_p\times S_{2p}$ we have
\ben(\sigma^{-1}\cdot \G)(x)_{\o\,
\o}=\sigma^{-1}(\G(x)_{{\sigma(\o)\, \sigma(\o)}})=
\sigma^{-1}(\G(x)_{{\o\, \o}}).
\end{equation*}
Hence
 $\G(x)_{\o\, \o}$ is $S_p\times S_{2p}$-invariant  proving
 \eqref{item:Im-G-is-W-invar}. The statement \eqref{item:diag-of-G(Gamma)}
 follows from \eqref{item:Im-G-is-W-invar} if we apply  a shift by $\ell\in
 \mP_0$ to an arbitrary $x\in \Y_p(\gl_2)$. The statement
 \eqref{item:diag-of-G(Gamma)-W-inw-defined-by-0} follows
 immediately from \eqref{item:Im-G-is-W-invar} and \eqref{item:diag-of-G(Gamma)}.
\end{proof}

 The composition $r_{\ell}=\Phi(\ell)\circ\G$ defines a
  representation
of $\Y_p(\gl_2)$.
By the construction, this representation provides a matrix realization
of the module $M(\ell)$; see Theorem~\ref{thm:basis}.

\bpr{}\label{prop:universal-G-is-faithful} The representation
$\G:$ $\Y_p(\gl_2){\longrightarrow}$ ${\M}_{{\mP}_0}(\bL)$ is
faithful. \epr

\begin{proof}
It is clear that
$$\Ker \G\subseteq \bigcap_{\ell\in
\mP_1}\Ker r_{\ell}.$$
Hence it is sufficient to prove that the intersection of the kernels is zero.
Let $\ell\in \mP_1$. Then $\Ker
r_{\ell}=\Ann M(\ell)$ by definition and so
 $\Ker r_{\ell}\subseteq
I_{\ell}$. However,
the intersection  $\bigcap_{\ell}I_{\ell}$ over all $\ell\in\mP_1$ is zero,
as was shown
in the proof of Proposition~\ref{cor:annih-all-generic=0}.
\end{proof}

\bco{}\label{cor:gamma-max-comm} \begin{enumerate}
\item\label{item:ga_is_maximal} ${\Gamma}$ is a maximal
commutative subalgebra in $\Y_p(\gl_2)$.
\item\label{item:diagonals_are_in_ga} If for $x\in \Y_p(\gl_2)$
the matrix $\G(x)$ is diagonal then $x\in {\Gamma}$.
\end{enumerate}\eco

\begin{proof} Consider an element $x\in$ $\Y_p(\gl_2)$ which commutes with
every $z\in$ ${\Gamma}$ and such that
$x\notin\Ga$. Suppose that there exist $\ell$, $\ell'\in
\mP_0$, $\, \ell\ne \ell'$ such that $\G(x)_{{\ell
}{\ell' }}\ne 0$. There exists $z\in {\Gamma}$ such that
$z(\ell)\ne$ $z(\ell')$ and so $\G(z)_{{ \ell }{\ell
}}\ne$ $\G(z)_{{\ell' }{\ell' }}$ by
Lemma~\ref{lem:W-invariants}\ts\eqref{item:diag-of-G(Gamma)}.
Then we
have
\ben
\G(xz)_{{\ell}{\ell'}}=\G(x)_{ { \ell} {\ell'}} \G(z)_{
{\ell'} {\ell'}}=\G(z x)_{ {\ell} {\ell'}}=\G(z)_{\ell {\ell
}} \G(x)_{ {\ell} {\ell'}}
\end{equation*}
which contradicts to the assumption.
Therefore $\G(x)$ is diagonal. To
prove the maximality of $\Gamma$ it is now sufficient to verify
part~\eqref{item:diagonals_are_in_ga} of the corollary.
By Lemma~\ref{lem:W-invariants}\ts\eqref{item:Im-G-is-W-invar},
we have $\G(x)_{\o\,\o}={f}/{g}\in K(\Ga)$ with relatively prime
$f,g\in{\Gamma}$. Suppose that $g\not\in
{\Bbbk}$. By Lemma~\ref{lem:W-invariants}\ts\eqref{item:Im-G-is-W-invar}
and \eqref{item:diag-of-G(Gamma)-W-inw-defined-by-0}, we derive that
$\G(x)\G(g)=$ $\G(f)$ and hence $xg=f$ by
Proposition~\ref{prop:universal-G-is-faithful}.
This shows that $x\in \Gamma$
due to Theorem \ref{thm:freedom}\ts\eqref{thm:Y-free-over-Gamma}.
\end{proof}

Denote by $X_0$ the column matrix defined by
\ben
X_0=\sum_{\delta\in \mP_0}\bL \ts e_{\delta,\o},
\end{equation*}
where
$\overline 0$ is the zero element of $\mP_0$.
Note that the $\bW$-action \eqref{waction}
induces an action of $S_p\times S_{2p}$ on the free
$\bL$-module $X_0$.

\bco{}\label{cor:im-pG-S-invariant} Let $p:$
${\M}_{\mP_0}({\bL}){\to}$ $X_0$ be the natural projection.
Then the composition $r=p\circ\G :\Y_p(\gl_2)\to X_0$ is injective.
Moreover, the map $p$ commutes with the action of
$S_p\times S_{2p}$ and, in particular, $r(\Y_p(\gl_2))$ is
$S_p\times S_{2p}$-invariant.
\eco

\begin{proof} Note that for any $x\in \Y_p(\gl_2)$ the matrix
 $\G(x)\in$ ${\M}_{\mP_0}(\bL)$ is determined completely by
its column $p(\G(x))$.
 Thus $r(x)=0$ implies
$\G(x)=0$ and hence $x=0$ since $\G$ is faithful. This proves that $r$ is injective.
The other statements follow immediately from the
definitions and Lemma~\ref{lem:W-invariants}\ts\eqref{item:Im-G-is-W-invar}.
\end{proof}

As in Section~\ref{sec:Some},
we identify the
$({\Gamma}{-}{\Gamma})$-bimodule structure on $\Y_p(\gl_2)$ with the
corresponding action of ${\Gamma}\otimes{\Gamma}$.
Using the embedding \eqref{imath}, we can regard
the elements of ${\Gamma}\otimes{\Gamma}$ as polynomials
in two families of variables $\mathsf b$ and $\mathsf b'$
which are $S_p\times S_{2p}$-invariant.

\ble{}\label{lem:acts-Gamma-otimes-Gamma} Suppose that $x\in
\Y_p(\gl_2)$, $f\in$ $\Gamma\otimes \Gamma$, and $\ell$, $\ell'\in$
$\mP_0$. Then

$$\G(f\cdot x)_{ {\ell} {\ell'}} \ = \
f(\mathsf b +\ell, \mathsf b+\ell') \G(x)_{ {\ell} {\ell'}}.$$
\ele

\begin{proof}
Let $f=\sum_{i}z_i\otimes z_i'\in \Gamma\otimes \Gamma$. Then
$\G(f\cdot x)=\sum_{i}\G(z_i)\G(x)\G(z_i')$ and hence,
by Lemma~\ref{lem:W-invariants}\ts\eqref{item:diag-of-G(Gamma)-W-inw-defined-by-0},
$$\G(f\cdot
x)_{ {\ell} {\ell'}}= \sum_{i}\G(z_i)_{\ell \, \ell}\G(x)_{\ell \,
\ell'}\G(z_i')_{\ell '\, \ell'}=\G(x)_{\ell \,
\ell'}\sum_{i}\G(z_i)_{\ell \, \ell}\G(z_i')_{\ell' \, \ell'}=$$
$$\G(x)_{\ell \, \ell'}\sum_{i}z_i(\mathsf
b+\ell)z_i'(\mathsf b+\ell')= \G(x)_{\ell \, \ell'}f(\mathsf b
+\ell, \mathsf b+\ell').$$
\end{proof}

\ble\label{lem:ann-Gamma-otimes-Gamma-poly}
  Let  $S\subseteq{\mL}$ be an $S_p\times
S_{2p}$-invariant set. Suppose that $z\in \Gamma$  and  $x\in \Y_p(\gl_2)$
is such that
$\G(x)_{\ell \, \ell'}=0$ for all $\ell, \ell'$,
$\ell-\ell'\not\in S$.  Then $F_{S,z}\cdot x=0$. \ele

\begin{proof} Let $F=F_{S,z}=\sum_i z^i\otimes a_i$ with
$a_i\in \L$. If $\ell-\ell'\in S$ then
\ben
\G(F\cdot x)_{ { \ell} {
\ell' }}=F(z(\mathsf b+\ell), \mathsf b+\ell')\ts\G(x)_{\ell \, \ell'}
\end{equation*}
by Proposition~\ref{prop:prop-sootn}(\ref{item:coeff-are-in-Gamma}) and
  Lemma~\ref{lem:acts-Gamma-otimes-Gamma}. Furthermore, observe that
$h(z, \mathsf b)=z\otimes 1 -
1\otimes z(\mathsf b+\ell-\ell')$ divides $F$ and that
$h(z(\mathsf b+\ell),\mathsf b+\ell')=0$.
Here we regard the result of the evaluation of the product of type
$z\ot z'(\mathsf b')$ at $\mathsf b$ as the polynomial
$z(\mathsf b)z'(\mathsf b')$.
This gives $F(z(\mathsf b+\ell), \mathsf b+\ell')=0$.
Hence, $\G(F\cdot x)=0$ implying
$F\cdot x=0\ts$ by Proposition~\ref{prop:universal-G-is-faithful}.
\end{proof}

 Let $S\subseteq$ $\mP_0$ be a  finite $S_p\times S_{2p}$-invariant
set. Define
$\Y^S$$=\{x\in
\Y_p(\gl_2)\,|\,{\mathscr L}_x\subseteq S\}$. Clearly $\Y^S$ is a
${\Gamma}$-subbimodule of $\Y_p(\gl_2)$.
 We have the following characterization of the bimodule $\Y^S$.

\ble{}\label{lem:YS-properties} Let $x\in \Y_p(\gl_2)$. Then
\begin{enumerate}

\item\label{item:YS-on-general} $x\in$ $\Y^S$ if and only if
the condition  $\G(x)_{ {\ell}, {\ell'}}\ne 0$, for some $\ell$,
$\ell'\in$ $\mP_0$,  implies that $\ell-\ell'\in S$.

\item\label{item:YS-by-polynomial}
 $F_{{\mathscr L}_x \setminus S,z}\cdot x\in$ $\Y^S$ for any $z\in \Gamma$.

\item\label{item:YS-is-fg-Gamma-mod} $\Y^S$ is a finitely
generated left (right) ${\Gamma}$-module and $\Y^S=$ $\mathbb
D(\Y^S)$.

\item\label{item:Y0=Gamma} $\Y^{\{0\}}=\Gamma$.
\end{enumerate}\ele

\begin{proof} The statement \eqref{item:YS-on-general}  follows from
the definition of $\Y^S$.
 Let $F=F_{{\mathscr L}_x \setminus S,z}$ and $y=F\cdot x$. To prove
\eqref{item:YS-by-polynomial} calculate the matrix element
$\G(y)_{\ell \ell'}$ provided that $\ell-\ell'\not\in S$.
By Lemma~\ref{lem:acts-Gamma-otimes-Gamma},
\ben
\G(y)_{\ell \ell'}=\G(F\cdot x)_{\ell \ell'}=F(z(\mathsf b+\ell),
\mathsf b+\ell')\G(x)_{\ell \ell'}.
\end{equation*}
If
$\ell-\ell'\not\in \mL_x$ then $\G(x)_{\ell \ell'}=0$ and hence
$\G(y)_{\ell \ell'}=0$. Suppose now that $\ell-\ell'\in {\mathscr L}_x
\setminus S$. Then
$$F(z(\mathsf b+\ell), \mathsf b+\ell')=
\prod_{\delta\in {\mathscr L}_x \setminus S}
(z(\mathsf b+\ell)-z(\mathsf b +\ell' +\delta))=0.$$
This proves \eqref{item:YS-by-polynomial}.

Let $x\in \mathbb D(\Y^S)$ and suppose that $z\in \Gamma$ is such that $z\neq
0$ and $z x\in \Y^S$. Since $\G(z x)_{\ell \ell'}=$ $z(\mathsf
b+\ell)\G(x)_{\ell \ell'}$ by Lemma~\ref{lem:acts-Gamma-otimes-Gamma},
we have $\G(z x)_{\ell \ell'}=0$ if and
only if $\G(x)_{\ell \ell'}=0$ implying that $x\in \Y^S$. Hence
$\Y^S=$ $\mathbb D(\Y^S)$.

Consider $r(\Y^S)$ as a ${\Gamma}$-submodule of $X_0$ where
$r:\Y_p(\gl_2)\to X_0$ is defined in Corollary~\ref{cor:im-pG-S-invariant}.
Then $r(\Y^S)$ belongs to the free
$\mathbb L$-submodule $\sum_{\ell\in
S}\mathbb L e_{{\ell}\ts {\overline 0} }$
of $X_0$ of finite rank. Hence $\mathbb L\cdot
r(\Y^S)$ is finitely generated $\mathbb L$-module. Without loss of
generality we can assume that this module is generated by the elements
$r(x_1), \ldots, r(x_s)\in r(\Y^S)$.
Since $\mathbb D(\Y^S)=\Y^S$  we have
\ben
\mathbb
D\Big( {\sum_{i=1}^s}{\Gamma} x_i\Big)\subseteq\Y^S.
\end{equation*}
Fix
$x\in \Y^S$. Then
\ben
r(x)= {\sum_{i=1}^s}t_i\ts
r(x_i), \qquad t_i\in \mathbb L.
\end{equation*}
Note that for any $y\in \Y^S$ and
any $\sigma\in S_p\times S_{2p}$ we have $\sigma\cdot r(y)=r(y)$
since $S$ is $S_p\times S_{2p}$-invariant.
Hence
\ben
p!\ts(2p)!\ts r(x)= \sum_{\sigma\in S_p\times S_{2p}}
\sigma \cdot r(x)= {\sum_{\sigma\in S_p\times
S_{2p}}}\ts  {\sum_{i=1}^s} (\sigma\cdot t_i)\ts\sigma\cdot
r(x_i)
\end{equation*}
which can be rewritten as
\beql{rx}
r(x)=
\frac{1}{p!\ts(2p)!}\ts{\sum_{i=1}^s}u_i\ts r(x_i),
\qquad
\text{where}\quad u_i=
 {\sum_{\sigma\in S_p\times S_{2p}}} \sigma\cdot
t_i.
\end{equation}
Since each $u_i$ is $S_p\times S_{2p}$-invariant, it
belongs to the field of fractions $K(\Gamma)$ for all $i=1,
\ldots, s$. Multiplying both parts of \eqref{rx} by the
common denominator of the $u_i$ we obtain from
Corollary~\ref{cor:im-pG-S-invariant} that
\ben
x\in\mathbb
D\Big( {\sum_{i=1}^s}{\Gamma} x_i\Big),\qquad\text{implying}\qquad
\mathbb
D\Big( {\sum_{i=1}^s}{\Gamma} x_i\Big)=\Y^S.
\end{equation*}
Due to
Corollary~\ref{cor:D(P)-is-fin-gen-over-center}, we can conclude that
$\Y^S$ is finitely generated over $\Gamma$. This proves
\eqref{item:YS-is-fg-Gamma-mod}. By the definition of $\Y^S$,
$x\in \Y^{\{0\}}$ if and only if $\G(x)$ is diagonal. Hence $x\in
\Gamma$ by Corollary~\ref{cor:gamma-max-comm}\ts\eqref{item:diagonals_are_in_ga}.
\end{proof}

\section{Category of Harish-Chandra modules}\label{sec:cat-of-HC-mod}
\setcounter{equation}{0}

We will show in this section that each character of $\Gamma$ extends to a finite number
of irreducible  Harish-Chandra modules over $\Y_p(\gl_2)$. This is
an analogue of the corresponding result in  the case of a Lie
algebra $\gl_n$ which was conjectured in \cite{dfo:gz} and proved
in \cite{o:fs}. In this section we use the techniques of
 \cite{dfo:hc} and \cite{o:fs}.

Since $\Gamma$ is a Harish-Chandra subalgebra of $\Y_p(\gl_2)$ we can apply
all the statements from Section~\ref{subsec:hca}.
Set ${\mathscr A}={\mathscr A}_{\Y_p(\gl_2),\Gamma}$. Then by
Proposition~\ref{prop:HC-mod-equiv-A-mod},
the categories ${\mathscr A}\dmo_d$
 and $\mathbb H(\Y_p(\gl_2),{\Gamma})$ are equivalent. Also
the full subcategory
 $\mathbb HW(\Y_p(\gl_2),{\Gamma})$ consisting of weight modules  is equivalent to
 the module category ${\mathscr A}_W\dmo$.
If $\ell\in \mL$  then the category
$R_{\ell}=\mathbb HW(\Y_p(\gl_2),{\Gamma},D(\ell))$ is equivalent to
the block ${\mathscr A}_W(D(\ell))\dmo$ of the category
${\mathscr A}_W\dmo$.

Let ${\mathbf m}$, ${\mathbf n}\in$ $\Sp{\Gamma}$, $\ell_{{\mathbf
m}}$,  $\ell_{{\mathbf n}}\in \mL$ are such that
$\imath^*(\ell_{{\mathbf m}})={\mathbf m}$ and
$\imath^*(\ell_{{\mathbf n}})={\mathbf n}$; see \eqref{imathstar}.
Set $$S({\mathbf
m},{\mathbf n}) \ = \{ \sigma_1\ell _{{\mathbf n}} -
 \sigma_2 \ell _{{\mathbf m}} \ |\  \sigma_1,\sigma_2\in  S_p\times S_{2p} \}
\cap \mP_0.$$
% tak kak peresekaem s L_0, to S_2p na samom dele ne deystvuet.
Consider
the following subset in ${\mL}$
$$\mP_2=\{\ell \in {\mathscr L} \ |\
\ell _{i}-\ell _{j}\not\in \mathbb Z \setminus
\{0 \},\qquad i,j=1,\ldots, p\}$$  and put
$\Omega=\imath^*(\mP_2)$. We shall also be using
the set ${\mathscr A}(\bm, \bn)$ introduced in \eqref{equ_2}.

\bpr{}\label{prop:main-YS-generates-A(m-n)+coroll}
\begin{enumerate}

\item\label{item:prop-main-Ys-Ys-generates} For any ${\mathbf m}$,
${\mathbf n}\in$ $\Sp{\Gamma}$ and any $m,n \geq 0$ we have
$$\Y_p(\gl_2)=\Y^S+{\mathbf n}^n \Y_p(\gl_2)+
\Y_p(\gl_2){\mathbf m}^m,$$ where $S=S({\mathbf m},{\mathbf n})$.

\item \label{item:prop-main-Ys-Gamma-big}
${\mathscr A}({\mathbf m},{\mathbf n})$ is
finitely generated as a left $\Gamma_{\mathbf n}$-module and as a right
${\Gamma}_{\mathbf m}$-module. In particular, the algebra
${\Gamma}$ is big in every ${\mathbf n}\in$ $\Ob{\mathscr A}$.

\item \label{item:prop-main-Ys-Yst0} If $S({\mathbf m},{\mathbf
n})=\{0\}$ then  ${\mathscr A}({\mathbf m},{\mathbf n})$  is
generated as a left $\Gamma_{\bn}$-module and as a right
$\Gamma_{\bm}$-module by the image
 of $1$ in ${\mathscr A}({\mathbf
m},{\mathbf n})$.

\item \label{item:prop-main-Ys-Ys-0} If $S({\mathbf m},{\mathbf
m})=\{0\}$ then $\,{\mathbf m}\in$ $\Omega$. Moreover, ${\mathscr
A}({\mathbf m},{\mathbf m})$ is a quotient algebra of $\Gamma_{\mathbf m}$ and
$\bm$ extends uniquely to an irreducible
$\Y_p(\gl_2)$-module.

\item \label{item:prop-main-Ys-generic-Gamma}If $\ell_{\mathbf
m}\in \mP_1$ then ${\mathscr A}({\mathbf m},{\mathbf m})=
\Gamma_{\mathbf m}$.

\item\label{item:gen-A(m-n)-free-rank-1} Let $\ell\in \mP_1$, $\bm
=\imath^*(\ell)$
 and $\bn=\imath^*(\ell+\delta_i)$, $i\in \{1, \ldots, p\}$.
 Then  ${\mathscr A}(\bm, \bn)$ is a free of rank
$1$ as a right  $\Gamma_{\bm}$-module and as a left $\Gamma_{\bn}$-module.
\end{enumerate}
\epr

\begin{proof} \eqref{item:prop-main-Ys-Ys-generates}
We shall show that for
 any $x\in \Y_p(\gl_2)$ and any
$k\geqslant 1$ there exists $x_k\in$ $\Y^S$ such that

\beql{equ_8}
x\in x_k+ \displaystyle {\sum_{i=0}^k}{\mathbf
n}^{k-i} x {\mathbf m}^i.
\end{equation}
The statement will follow if we choose $k=m+n+1$. We will use
induction on $k$. Suppose that $k=1$.
If $\mL_x\subseteq S$ then $x\in \Y^S$ and there
is nothing to prove. Furthermore, by the definition of the set $S$
for any $\ell \in$ ${\mathscr L}_x \setminus S$ the $S_p\times
S_{2p}$-orbits of $\ell_{\bn}$ and  $\ell_{\bm}+\ell$ are
disjoint. Hence there exists $z\in{\Gamma}$  such that $z(\ell
_{{\mathbf n}})\ne$ $z(\ell _{{\mathbf m}}+\ell )$ for any $\ell
\in$ ${\mathscr L}_x \setminus S$.
 Let $F=F_{{\mathscr L}_x
\setminus S,z}$. Then
\ben
F(z(\ell_{\mathbf n}), \ell_{\mathbf m})=
\prod_{\ell \in{\mathscr L}_x \setminus S} (z(\ell _{{\mathbf
n}})-z(\ell _{{\mathbf m}}+\ell ))\neq 0.
\end{equation*}
We can assume that $F(z(\ell_{\mathbf n}),
\ell_{\mathbf m})=1$. Hence we obtain that
 $F=1+u$ where $u\in{\mathbf n}\otimes\Gamma+$
$\Gamma\otimes{\mathbf m}$. It follows from
Lemma~\ref{lem:YS-properties}\ts\eqref{item:YS-by-polynomial},  that
$x_1=F\cdot x$ belongs to $\Y^S$. Hence we have
\ben
x_1=(1+u)\cdot
x\in x+{\mathbf n}\ts x\ts \Gamma + \Gamma\ts  x\ts {\mathbf m}\qquad
\text{and thus}\qquad
x\in
x_1+{\mathbf n}\ts  x\ts \Gamma+\Gamma\ts  x\ts {\mathbf m}.
\end{equation*}
This proves the assertion in the case $k=1$.
Assume that \eqref{equ_8} holds for some $k\geq
1$. Then
$$x\in x_k + \displaystyle {\sum_{i=0}^k} {\mathbf n}^{k-i} (x_k +
\displaystyle {\sum_{j=0}^k} {\mathbf n}^{k-j} x{\mathbf m}^j)
{\mathbf m}^i \subseteq x_k + \displaystyle {\sum_{i=0}^k} {\mathbf
n}^{k-i} x_k {\mathbf m}^i+ \displaystyle {\sum_{i=0}^{k+1}}
{\mathbf n}^{k+1-i} x {\mathbf m}^i.$$ Since $\Y^S$ is a
$\Gamma$-bimodule we conclude that $x_k + \displaystyle
{\sum_{i=0}^k} {\mathbf n}^{k-i} x_k {\mathbf m}^i\subseteq \Y^S$
which implies the statement
\eqref{item:prop-main-Ys-Ys-generates}. In particular,
we have proved that
\beql{xkxk-1}
x_{k+1}-x_k\in \sum_{i=0}^k {\mathbf n}^{k-i} \Y^S {\mathbf m}^i.
\end{equation}

\eqref{item:prop-main-Ys-Gamma-big} We prove the statement for the
case of left module, the case of right module can be treated
analogously. By \eqref{item:prop-main-Ys-Ys-generates} the image
$\overline x$ of every $x\in \Y^S$ in ${\mathscr A}({\mathbf
n},{\mathbf m})$ is the limit of
 the  sequence $(\overline{x}_k)_{k\geqslant
1}$, $x_k\in$ $\Y^S$. Let $y_1, \ldots, y_m$ be a finite system of
generators of $\Y^S$ as a left ${\Gamma}$-module which exists by
Lemma~\ref{lem:YS-properties}\ts\eqref{item:YS-is-fg-Gamma-mod}.
Then for every
$N>1$ and every $i=1, \ldots, m$
 there exists  $N_i$ such that
$$y_i{\mathbf m}^N\subseteq \sum_{j=1}^m {\mathbf n}^{N_i}y_j.$$
Since $\Gamma$ is noetherian we have that $\bigcap_{k}\bn^k\Y^S=0$
and hence there exists the maximum value
 $d_N$ such that
$$y_i{\mathbf m}^N\subseteq \sum_{j=1}^m {\mathbf n}^{d_N}y_j$$
for all $i=1, \ldots, m$. Moreover,  $d_N\rightarrow
\infty$ while $N\rightarrow \infty$ since $\Y^S$ is a finitely generated right
$\Gamma$-module and
$\bigcap_{k}\Y^S\bm^k=0$.
  By \eqref{xkxk-1}, $x_{k+1} - x_k\in
 {\mathbf n}^{R_k} \Y^S$ where $R_k=\min\{[k/2],
d_{[k/2]}\}$.     We have
$$\overline x=\overline x_1 + \sum_{k=1}^{\infty}\overline{(x_{k+1}-x_k)}$$
 and thus
\ben
\overline x\in \sum_{k=1}^{\infty}\overline{
{\mathbf n}^{R_k}\Y^S}\subseteq \sum_{l=1}^m
 \Gamma_{\mathbf n}\ts\overline {y_l}.
\end{equation*}
 Note that the first sum is well defined since
$R_k\rightarrow \infty$ when $k\rightarrow \infty$.
  We conclude that  ${\mathscr A}({\mathbf
n},{\mathbf m})$ is finitely generated  as a left
$\Gamma_{\bn}$-module.  This completes the proof of
\eqref{item:prop-main-Ys-Gamma-big}.

\eqref{item:prop-main-Ys-Yst0}
By Lemma~\ref{lem:YS-properties}\ts\eqref{item:Y0=Gamma},
 $\Y^{\{\,0\}}=\Gamma$. Hence
 $\mathcal A(\bm,\bn)$ is generated (both as a left and
as a right module) by the image of $1\in\Gamma$ by
(\ref{item:prop-main-Ys-Ys-generates}).

\eqref{item:prop-main-Ys-Ys-0} By \eqref{item:prop-main-Ys-Yst0},
 $\mathcal A(\bm,\bm)$ is $1$-generated as a
left $\Gamma_{\bm}$-module.  Then the $\Bbbk$-algebra homomorphism
\ben
{\hat \imath}_{{\mathbf m}}:{\Gamma}_{{\mathbf
m}}{\to}\mathcal A(\bm,\bm),\qquad z\mapsto z\cdot
\1_{\bm},
\end{equation*}
where $\1_{\bm}$ is a unit morphism, is an epimorphism
which shows that $ A(\bm,\bm)$ is a quotient algebra of
${\Gamma}_{\mathbf m}$. The uniqueness of the extension follows from
the uniqueness of the simple $A(\bm,\bm)$-module and
\cite[Theorem~18]{dfo:hc}.

\eqref{item:prop-main-Ys-generic-Gamma} Let $\ell =\ell _{{\mathbf
m}}\in \mP_1$. Then $S(\bm,\bm)=\emptyset$ and
 $A(\bm,\bm)\simeq \Gamma_{\bm}/J_{\bm}$ by (\ref{item:prop-main-Ys-Ys-0})
and thus $J_{\bm}$ acts trivially on $M(\bm)$ in any Harish-Chandra module
$M$.   Since $\ell\in\mP_1$ then for any $k>0$  there exists a
canonical projection $\tilde{\pi}_{k}:$ $\mathbb L{\longrightarrow}$
$\mathbb L/ {(\ell)}^k$, where $(\ell)^k=\ell^k \mathbb L$.
It induces a homomorphism of the
matrix algebras $\pi_k:{\M}_{\mP_0}(\bL)\longrightarrow
{\M}_{\mP_0}(\bL/ {(\ell)}^k)$ and defines a Harish-Chandra module
by the following composition
$$\Y_p(\gl_2){\xar{\G}} {\M}_{\mP_0}(\bL)
\xar{\pi_k}{\M}_{\mP_0}(\bL/ {(\ell)}^k).$$ For any nonzero
 $x\in \Gamma$
there exists $k>0$ such that $x\not\in (\ell)^{k}$ and hence
 $\pi_k
(\G(x)_{\o,\o})= x + (\ell)^k\ne0$. Therefore, there exists
a Harish-Chandra module $M$ where $x$ acts nontrivially on $M(\bm)$
implying that $J_{\bm}=0$.  This completes the proof.

\eqref{item:gen-A(m-n)-free-rank-1} The proof is analogous to the
proof of \eqref{item:prop-main-Ys-generic-Gamma}. Let $z\in
\Gamma, $ $z\neq 0$. Suppose ${\mathscr A}(\bm, \bn)z=0$. Then by
the construction of the equivalence  $\mathbb F:{\mathscr
A}\dmo_d$ ${\longrightarrow}$ $\mathbb H(U,{\Gamma})$ for any
Harish-Chandra module $M$ and any $x\in \Y_p(\gl_2)$ the linear
operator $xz$ on $M$ induces the zero map between $M(\bm)$ and
$M(\bn)$. It is sufficient to construct a Harish-Chandra module where
this is failed. For $k\geq 1$ consider as in (\ref{item:prop-main-Ys-generic-Gamma})
the composition $\pi_k\circ \G:\Y_p(\gl_2)\to
\M_{\mP_0}(\bL/(\ell)^k)$. It defines a Harish-Chandra module
structure on a free $\bL/(\ell)^k$-module
\ben
\overline{X}=\sum_{\delta\in
\mP_0}\bL/(\ell)^k e_{\delta,\o}.
\end{equation*}
Consider $x\in \Y_p(\gl_2)$
such that $\G(x)_{\delta_i \overline{0}}\neq 0$ for some $i$. Then
\ben
\G(xz)_{\delta_i \overline{0}}=\G(x)_{\delta_i
\overline{0}}\G(z)_{\overline{0} \overline{0}}=\G(x)_{\delta_i
\overline{0}}z\neq 0.
\end{equation*}
Choose  $k$ such that $\G(xz)_{\delta_i
\overline{0}}\not\in (\ell)^k$. Hence
 $(\pi_k\cdot \G)(xz)_{\delta_i, \overline{0}}\neq 0$ and
 the linear operator $xz$  induces
a non-zero map between $\overline{X}(\bm)=\bL/(\ell)^k$ and
$\overline{X}(\bn)=\bL/(\ell+\delta_i)^k$. The contradiction
shows that ${\mathscr A}(\bm, \bn)z\neq 0$. The case $z{\mathscr
A}(\bm, \bn)=0$ is treated in a similar manner.
\end{proof}

Now we are in a position to state the main result of this
section which follows immediately from Lemma \ref{lem:Gamma-big}
and
Proposition~\ref{prop:main-YS-generates-A(m-n)+coroll}\ts
\eqref{item:prop-main-Ys-Gamma-big}.

\bth{}\label{thm:main} Let ${\mathbf m}\in$ $\Sp{\Gamma}$. Then
the left ideal $\Y_p(\gl_2){\mathbf m}$  is contained in finitely
many maximal left ideals of $\Y_p(\gl_2)$. In particular,
${\mathbf m}$ extends to a finitely many (up to an
isomorphism) irreducible $\Y_p(\gl_2)$-modules and for each such
module $M$,
   $\dim M({\mathbf n})<$ $\infty$ for all ${\mathbf n}\in$ $\Sp{\Gamma}$.
 \eth

\section{Category of generic Harish-Chandra modules}\label{sec:Generic}
\setcounter{equation}{0}

In this section
we study a full
subcategory of generic modules
in {$ \mathbb HW(\Y_p(\gl_2),{\Gamma})$}.
We give a complete description of irreducible modules
and indecomposable modules in tame
blocks of this category.

\ble\label{lem:isoob} Let $\ell\in \mP_1$, $\ell=(\beta, \gamma)$,
$\bm=\imath^*(\ell)\in \Sp \Gamma$, $\bn=\imath^*(\ell+\delta_i)$,
$i\in \{1, \ldots, p\}$. If $\beta_i\not\in \{\gamma_1, \ldots,
\gamma_{2p}\}$ then the objects of ${\mathscr A}$ represented by
$\bm$ and $\bn$ are isomorphic. \ele

\begin{proof}
Choose $z_1, z_2\in \Gamma$ such that
\ben
z_1(\ell+\delta_j)=\delta_{ij},\qquad
z_2(\ell+\delta_i-\delta_j)=\delta_{ij},
\qquad j=1, \ldots, p.
\end{equation*}
Set $z=z_2\ts t_{12}^{(1)}z_1\ts t_{21}^{(1)}$.  Then $G(z)$ is
diagonal by Lemma \ref{lem:acts-Gamma-otimes-Gamma} and hence
$z\in \Gamma$ by Corollary \ref{cor:gamma-max-comm}\ts\eqref{item:diagonals_are_in_ga}.
We will show that the image of
$z$ in $\Gamma_{\bm}$ is invertible. Clearly, this is equivalent
to the fact that $z(\ell)\neq 0$. Formulas \eqref{T22act}--\eqref{t11u}
imply that
$z(\ell)=\gamma(-\beta_i)\neq 0$ by assumption.   Denote by
$T_1$ (respectively $T_2$) the generator of
$\Gamma_{\bm}-\Gamma_{\bn}$(respectively, $\Gamma_{\bn}-\Gamma_{\bm}$)-bimodule
${\mathscr A}(\bm, \bn)$ (respectively,
${\mathscr A}(\bn, \bm)$); see
Proposition~\ref{prop:main-YS-generates-A(m-n)+coroll}\ts
\eqref{item:gen-A(m-n)-free-rank-1}.
Then
\ben
z_2\ts t_{12}^{(1)}=z_{\bm
}T_2,\qquad  z_1\ts t_{21}^{(1)}=T_1z'_{\bm}
\end{equation*}
for some $z_{\bm},z'_{\bm}\in
\Gamma_\bm$ and hence $z=z_{\bm}T_2 T_1z'_{\bm}$. Since $z(\ell)\ne 0$ it
follows that $z'_{\bm}(\ell)\ne 0,z_{\bm}(\ell)\ne 0$ and so $T_2
T_1= z_{\bm}^{-1}z (z'_{\bm})^{-1}$ is invertible in
$\Gamma_{\bm}$.
A similar argument shows that
$T_1T_2$ is  invertible in $\Gamma_{\bn}$. Therefore the objects
$\bm$ and $\bn$ are isomorphic.
\end{proof}

\bco\label{cor:ext-strongly-gen} Let  $\ell\in \mP_1$,
$\ell=(\beta, \gamma)$, $\beta_i-\gamma_j\not\in \Z$ for all $i,j$.
Then the category $\mathbb
H(\Y_p(\gl_2),{\Gamma}, D(\ell))$ is hereditary. Moreover,
$$\dim \Ext^1_{\mathbb H(\Y_p(\gl_2),{\Gamma},D(\ell))}(L(\ell), L(\ell))=3p.$$
\eco

\begin{proof}
  Let
$\bm=\imath^*(\ell)\in \Sp \Gamma$.
By Lemma \ref{lem:isoob} and our assumptions all objects of the
category ${\mathscr A}(D(\ell))$ are isomorphic and hence the
category ${\mathscr A}(D(\ell))\dmo_d$ is equivalent to the
category of finite-dimensional modules over $\Gamma_{\bm}$.
Applying
 Proposition~\ref{prop:HC-mod-equiv-A-mod} we conclude that the category
 $\mathbb H(\Y_p(\gl_2),{\Gamma},
D(\ell))$ is hereditary. Since $\Gamma_{\bm}$ is an algebra of
power series in $3p$ variables the statement about $\dim \Ext^1$
follows.
\end{proof}

\subsection{Category of generic weight modules}\label{sec:Genericweight}
Let us fix
 \ben
\ell\in \mP_1,\qquad \bm=\imath^*(\ell),\qquad
\bn=\imath^*(\ell+\delta_i)\in \Sp \Gamma, \qquad i\in \{1, \ldots,
p\}.
\end{equation*}
Then ${\mathscr A}_W(\bm, \bm)\simeq
\Gamma_{\bm}/\Gamma_{\bm}\bm\simeq \Bbbk$ by
Proposition~\ref{prop:main-YS-generates-A(m-n)+coroll}\ts
\eqref{item:prop-main-Ys-generic-Gamma} and so, $\dim \mathscr
A_W(\bm, \bn)=1$ by Proposition~\ref{prop:main-YS-generates-A(m-n)+coroll}\ts
\eqref{item:gen-A(m-n)-free-rank-1}. We will give  a direct
construction of the category ${\mathscr A_W(D(\ell))}$.

We shall keep the notation
\ben
\ell=(\beta, \gamma), \qquad \be=(\beta_1, \ldots, \beta_p)\in
\Bbbk^p, \qquad \ga=(\ga_1, \ldots, \ga_{2p})\in \Bbbk^{2p}.
\end{equation*}
Since $\ell\in \mP_1$ then $\beta_i-\beta_j \notin \Z$
for $i\neq j$.
 Consider the following category $K_{\ell}$: ${\rm
Ob}\ts(K_{\ell})= \Z^p$ and the morphisms are generated by
\ben
f_i(k): (k)\mapsto (k+\delta_i) \qquad\text{and}\qquad
e_i(k): (k) \mapsto (k-\delta_i),
\end{equation*}
where
$i=1, \ldots, p$ and
$(k)=(k_1, \ldots, k_p)\in
\Z^p$ with the following relations:
\ben
\bal
f_j(k+\delta_i)\ts f_i(k)&= f_i(k+\delta_j)\ts f_j(k),\\
e_j(k-\delta_i)\ts e_i(k)&= e_i(k-\delta_j)\ts e_j(k),\\
e_i(k+\delta_j)\ts f_j(k)&=
f_j(k-\delta_i)\ts e_i(k)\qquad\text{for}\quad i\neq j,\\
e_i(k+\delta_i)\ts f_i(k)&=-\ga(-\beta_i-k_i)\ts 1_{(k)},\\
f_i(k-\delta_i)\ts e_i(k)&=-\ga(-\beta_i-k_i+1)\ts 1_{(k)}.
\eal
\end{equation*}
It follows immediately from Lemmas~\ref{lem:wgen} and
\ref{lem:comr} that any module in the category $R_{\ell}$ can
 be naturally viewed as a module over the category $K_{\ell}$ which
 defines a functor $F: R_{\ell}\rightarrow K_{\ell}\dmo$.
Consider the cyclic subalgebra $C_{\ell}(a)= {\rm
Hom}^{}_{K_{\ell}}(a, a)$ for any $a\in \Z^p$. Clearly,
$C_{\ell}(a)\simeq \Bbbk$ for any $a\in \Z^p$ due to the defining
relations of $K_{\ell}$. For any $a=(k_1, \ldots, k_p)\in \Z^p$ we
can construct a universal module $M(\ell, a)\in K_{\ell}\dmo$.
Consider $\Bbbk$ as a $C_{\ell}(a)$-module with
\ben
\bal
e_i(k+\delta_i)\ts f_i(k)\ts 1&=-\ga(-\beta_i-k_i),\\
f_i(k-\delta_i)\ts e_i(k)\ts 1&=-\ga(-\beta_i-k_i+1).
\eal
\end{equation*}
Let $A_{\ell,a}$ be an algebra of paths in $K_{\ell}$ originating
in $a$. Now construct a $\Z^p$-graded $K_{\ell}$-module
$$M(\ell, a)=A_{\ell, a}\otimes_{C_{\ell}(a)} \Bbbk.$$
Clearly, all graded components of $M(\ell, a)$ are $1$-dimensional
and $M(\ell, a)_a=1_a\otimes \Bbbk$. A module $M(\ell, a)$
contains a unique maximal $\Z^p$-graded submodule which intersects
$M(\ell, a)_a$ trivially and hence has a unique irreducible
quotient $L(\ell, a)$ with $L(\ell, a)_a\simeq \Bbbk$ and $\dim
L(\ell, a)_b\leq 1$ for all $b\in \Z^p$. If $V$ is another
irreducible $K_{\ell}$-module with $V_a\neq 0$ then there exists a
non-trivial $C_{\ell}(a)$-homomorphism from $\Bbbk$ to $V_a$ which
can be extended to an epimorphism from $M(\ell, a)$ to $V$. Since
$V$ is irreducible we conclude that $V\simeq L(\ell, a)$.

Obviously, we can view $M(\ell)$ as a module over the category
$K_{\ell}$ with a natural action of the morphisms of $K_{\ell}$
and $F(M(\ell))=M(\ell, \be)$. Thus a $K_{\ell}$-module $M(\ell,
\be)$ can be extended to a $\Y_p(\gl_2)$-module $M(\ell)$.
Moreover, the functor $F$ preserves the submodule structure of
$M(\ell)$. In particular, $F(L(\ell))=L(\ell, \be)$.

\bpr\label{prop:p6} If $\ell\in \mP_1$ then the categories
$K_{\ell}\dmo$ and $R_{\ell}$ are equivalent. \epr

\begin{proof} Let $\ell=(\beta, \ga)$.
We already have the functor $F: R_{\ell}\rightarrow K_{\ell}\dmo$.
Suppose that $V\in K_{\ell}\dmo$. We want to show
that $V$ can be extended to a $\Y_p(\gl_2)$-module. Fix $v\in
V_{(k)}\setminus \{0\}$. Let $W\subseteq V$ be a submodule
generated by $v$. Then $W_{(k)}=\Bbbk v$ and there is an
epimorphism from $M(\ell, a)$ to $W$, where $a=(k_1, \ldots,
k_p)$, which maps $1_a\otimes 1$ to $v$. Since  $F(M(\ell'))=
M(\ell, a)$, where $\ell'=(\beta + a, \ga)$, then $W$ can be
extended to a corresponding quotient of $M(\ell')$. Since $v$ was
an arbitrary element of $V$ we conclude that $V$ can be extended
to a $\Y_p(\gl_2)$-module and will denote that module by $G(V)$.
Clearly, $G$ defines a functor from $K_{\ell}\dmo$ to
$R_{\ell}$ (action on morphisms is obvious). One can easily see
that the functors $F$ and $G$ define an equivalence between the
categories $K_{\ell}\dmo$ and $R_{\ell}$.
\end{proof}

\subsection{Support of irreducible generic weight  modules}\label{sec:supp}

To complete the classification of irreducible modules we have to
know when two irreducible modules $L(\ell)$ and $L(\ell')$ are
isomorphic. For that we need to describe the support $\Supp
L(\ell)$.

We shall say that the weight subspaces $M(\ell)_{\psi}$ and
$M(\ell)_{\psi+\delta_i}$ are {\it strongly isomorphic\/} if
$\ga(-\psi_i)\neq 0$ where $\psi=(\psi_1, \ldots, \psi_p)$. This
implies
\ben
f_i(\psi_1, \ldots, \psi_p)\ts
M(\ell)_{\psi}\neq 0\quad\text{and}\quad e_i(\psi_1,
\ldots,\psi_i+1, \ldots, \psi_p)\ts M(\ell)_{\psi+\delta_i}\neq 0.
\end{equation*}

The  statement below follows immediately from the relations in
$K_{\ell}$ (cf. also Corollary~\ref{cor:gasupp}).

\ble\label{lem:stisom} If $M(\ell)_{\psi}$ and
$M(\ell)_{\psi+\delta_i}$ are strongly isomorphic, then
$M(\ell)_{\psi\pm \delta_j}$ and $M(\ell)_{\psi+\delta_i\pm
\delta_j}$ are strongly isomorphic for all $i,j=1, \ldots, p$ such that
$i\neq j$. Moreover, if
\ben f_i(\psi_1, \ldots, \psi_p)\ts
M(\ell)_{\psi}=0 \qquad\text{or}\qquad e_i(\psi_1, \ldots,
\psi_p)\ts M(\ell)_{\psi}=0
\end{equation*}
then
\ben
\bal
f_i(\psi_1, \ldots,\psi_j\pm 1, \ldots,
\psi_p)\ts M(\ell)_{\psi\pm \delta_j}&= 0
\qquad\text{or}\\
e_i(\psi_1, \ldots, \psi_j + 1, \ldots,  \psi_p)\ts
M(\ell)_{\psi\pm \delta_j}&=0,
\eal
\end{equation*}
respectively, for all $j\neq i$. \ele

Let $a_i, a'_i\in \Z\cup \{\pm \infty\}, a_i\leq a'_i$, $i\in \{1,
\ldots, p\}$. Denote
\ben
P(a_1, \ldots, a_p, a'_1,
\ldots, a'_p)= \{(x_1, \ldots, x_p)\in \Z^p\ |\  a_i\leq x_i\leq
a'_i,\ i=1, \ldots, p\},
\end{equation*}
a parallelepiped in $\Z^p$. Note that some faces of the
parallelepiped can be infinite in some directions. In particular,
in the case $a_i=-\infty$, $a'_i=\infty$ for all $i$, the
parallelepiped coincides with $\Z^p$.

\bth\label{thm:suppar} For any irreducible weight module $L(\ell)$
over $\Y_p(\gl_2)$ there exist elements $a_i, a'_i\in \Z \cup \{\pm
\infty\}$, $a_i\leq a'_i$, $i\in \{1, \ldots, p\}$ such that
\ben
\Supp L(\ell)=P(a_1, \ldots, a_p, a'_1, \ldots,
a'_p).
\end{equation*}
\eth

\begin{proof} Let $\ell=(\be,\ga)\in \mP_1$.
Fix $i\in \{1, \ldots, p\}$. If $\ga(-\beta_i+k)\neq 0$ for all
$k\in \Z$ then
\ben
(k_1, \ldots,  k_i+m,  \ldots,  k_p)\in
\Supp L(\ell)
\end{equation*}
as soon as $(k_1, \ldots, k_p)\in  \Supp L(\ell)$. This follows
immediately from Lemma~\ref{lem:stisom}. In this case we set
$a_i=-\infty$ and $a'_i=\infty$. Now let $\ga(-\beta_i+k)=0$ for
some $k\in \Z$. Let $m\geq 0$ be the smallest integer (if it
exists) such that $\ga(-\beta_i-m)=0$ and let $n\leq 0$ be the
largest integer (if it exists) such that $\ga(-\beta_i -n+1)=0$.
It follows from Lemma~\ref{lem:stisom} that \ben \Supp L(\ell)\cap
\{\be + k\delta_i\ |\ k\in \Z\}= \{\be+n\delta_i, \ldots, \be,
\ldots, \be+m\delta_i\}.
\end{equation*}
If $\be+s\delta_j\in \Supp L(\ell)$, $j\neq i$ then
\ben
\Supp L(\ell)\cap \{\be + s \delta_j + k\delta_i\ |\ k\in
\Z\}=\{\be+s\delta_j + n\delta_i, \ldots, \be+s\delta_j, \ldots,
\be+ s\delta_j+ m\delta_i\}.
\end{equation*}
In this case we set $a_i=\beta_i+n$ and $a'_i=\beta_i+m$. The
statement of the theorem now follows.
\end{proof}

\subsection{Indecomposable generic weight modules}\label{sec:ind}

Fix $\ell=(\be,\ga)\in \mP_1$. A full subcategory ${\mathcal
S}\subseteq K_{\ell}$ is called a {\it skeleton} of $K_{\ell}$
provided the objects of ${\mathcal S}$ are pairwise non-isomorphic
and any object of $K_{\ell}$ is isomorphic to some object of
${\mathcal S}$. In this case the categories of $K_{\ell}\dmo$ and
${\mathcal S}\dmo$ are equivalent.

For each $i\in \{1, \ldots, p\}$ consider a set $I_i=\{k\in \Z\ |\
\ga(-\beta_i-k)=0\}$. Define a category $S_{\ell}$ as a
$\Bbbk$-category with the set of objects
\beql{sob}
S_0=\{0,
\ldots, |I_1|\}\times \cdots \times \{0, \ldots, |I_p|\}
\non
\end{equation}
and with morphisms generated by

\beql{morso}
\bal
r_{(i_1, \ldots, i_p)}^k: (i_1, \ldots, i_p)
&\mapsto (i_1, \ldots, i_k+1, \ldots, i_p),\\
s_{(j_1, \ldots, j_p)}^k: (j_1, \ldots, j_p) &\mapsto (j_1,
\ldots, j_k-1, \ldots, j_p),
\eal
\non
\end{equation}
where $k\in \{1, \ldots, p\}$ is such that $I_k\neq \emptyset$,
$i_k< |I_k|$, $j_k>0$, subject to the relations:
\beql{relso}
s^k_{(i_1, \ldots, i_k+1, \ldots, i_p)} r^k_{(i_1, \ldots,
i_p)}=r^k_{(i_1, \ldots, i_p)} s^k_{(i_1, \ldots, i_k+1, \ldots,
i_p)}=0
\non
\end{equation}
and
$$
x^k_{(a_1, \ldots, a_p)}y^r_{(e_1, \ldots, e_p)}= y^r_{(c_1,
\ldots, c_p)}x^k_{(e_1, \ldots, e_p)}
$$
for all $k\neq r$ and all possible $x,y\in \{r, s\}$ and all $a_i, e_i,
c_i$, with $1\leq i\leq p$ for which this equality makes sense.

It follows from the construction that $S_{\ell}$ is the skeleton
of the category $K_{\ell}$. Note that the corresponding algebra is
finite-dimensional. In particular, $S_{\ell}$ is semisimple when
$I_k=\emptyset$ for all $1\leq k\leq p$, i.e. when $\ga(-\beta_k+
r)\neq 0$ for all $k\in \Z$ and all $i=1, \ldots, p$. Hence it is
sufficient to describe all indecomposable modules over $S_{\ell}$.

Fix $a\in S_0$ and define a simple $S_{\ell}$-module $S_a$ such
that $S_a(b)=\delta_{a,b}\Bbbk$ for all $b\in S_0$ and all
morphisms are trivial. Since $S_{\ell}$ defines a
finite-dimensional algebra we have the following

\bpr\label{prop:simple} Any simple module over $S_{\ell}$ is
isomorphic to $S_a$ for some $a\in S_0$. \epr

This is another confirmation of the fact that all weight spaces in
any irreducible generic weight $\Y_p(\gl_2)$-module are
$1$-dimensional. But this need not be the case for
indecomposable modules. We restrict ourselves to a full
subcategory $R_{\ell}^f\subseteq R_{\ell}$ which consists of
weight modules $V$ with $\dim V_{\psi}<\infty$ for all $\psi\in
\Supp V$. We will establish the representation type of the
category $R_{\ell}^f$ (finite, tame or wild). For the necessary
definitions we refer the reader to  \cite{d:tw}.

To establish the representation type of the category $R^f_{\ell}$
it is sufficient  to consider the category $S_{\ell}\dmo^f$ of
modules over the category $S_{\ell}$ with finite-dimensional
weight spaces. Denote $X_{\ell}=\{k\in \{1, \ldots, p\}\ |\
I_k\neq \emptyset\}$.

\subsubsection{Indecomposable modules in the case $|X_{\ell}|=1$}

In this section we describe all indecomposable modules over
$S_{\ell}$ in the case $|X_{\ell}|=1$. Let $X_{\ell}=\{i\}$ and
let $|I_i|=r>0$. In this case the category $S_{\ell}$ has the
following quiver ${\bf A}$ with the relations:

\begin{picture}(0.00,40.00)
%\put(00.00,17.00){${\bf A}:$}
\put(87.00,27.00){$a_1$}
\put(87.00,6.50){$b_1$}
\put(69.00,27.00){$1$}
\put(101.00,27.00){$2$}
\put(140.00, 27.00){$a_r$}
\put(140.00, 6.50){$b_r$}
\put(160.00,27.00){$r+1$}
\put(160.00,17.00){$\circ$}
\put(240.00,17.00){$a_i\ts b_i=b_i\ts a_i=0$}
\put(125.00,27.00){$r$}
\put(125.00,17.00){$\circ$}
\put(69.50,17.00){$\circ$}
\put(101.00,17.00){$\circ\ldots$}
\put(74.00,22.00){\vector(1,0){27.00}}
\put(101.50,16.50){\vector(-1,0){27.00}}
\put(130.00,22.00){\vector(1,0){27.00}}
\put(158.00,16.50){\vector(-1,0){27.00}}
\end{picture}

We denote by $S_i$, $i\in\{1,\ldots, r+1\}$, the simple module
corresponding to the point $i$. These modules correspond to all
irreducible modules in $R^f_{\ell}$ by
Proposition~\ref{prop:simple}. Now describe the remaining
indecomposable modules for the quiver above. Fix integers $1\leq
k_1<k_2\leq r+1$ and a function
\ben
\xi_{k_1,k_2}:\{k_1, k_1+1,
\ldots, k_2\} \rightarrow \{0,1\}.
\end{equation*}
Define the module $M=M(k_1,k_2,
\xi_{k_1,k_2})$ as follows: $M(i)=\Bbbk e_i$, $k_1\leq i\leq k_2$ and
$M(i)=0$ otherwise. Furthermore,
\begin{align}
a_ie_i&=e_{i+1},\qquad b_ie_{i+1}=0\qquad \text{if}\qquad
\xi_{k_1,k_2}(i)=1 \non\\
\intertext{and}
a_ie_i&=0, \qquad b_ie_{i+1}=e_i \qquad \text{if}\qquad
\xi_{k_1,k_2}(i)=0,     \non
\end{align}
for all $1\leq i<k_2$.
The proof of the following proposition is standard;
see e.g. \cite{gr:rf}.

\bpr\label{prop:indec}
The modules $S_i$ for $1\leq i\leq r+1$ and
$M(k_1,k_2,\xi_{k_1,k_2})$ with $1\leq k_1<k_2\leq r+1$ and
\beql{indec}
\xi_{k_1,k_2}:\{k_1, k_1+1, \ldots, k_2\}\rightarrow \{0,1\},
\non
\end{equation}
exhaust all non-isomorphic
indecomposable modules for ${\bf A}$.
\epr

\subsubsection{Indecomposable modules in the case $|X|_{\ell}=2$}

In this section we describe the indecomposable modules for
$S_{\ell}$ when $|X|_{\ell}=2$ and $|I_k|=1$ for each $k\in
X_{\ell}$. Then $S_{\ell}$ is isomorphic to the following category
${\bf B}$ considered in \cite{bb:ir}.

\begin{picture}(0.00,85.00)
\put(00.00,35.00){${\bf B}:$} \put(89.00,35.00){$a_{2}$}
\put(109.00,35.00){$b_{2}$} \put(48.50,35.00){$b_{0}$}
\put(67.00,35.00){$a_{0}$} \put(80.00,68.00){$a_{1}$}
\put(80.00,47.50){$b_{1}$} \put(80.00,26.50){$b_{3}$}
\put(80.00,7.50){$a_{3}$} \put(49.00,57.00){$1$}
\put(110.00,57.00){$2$} \put(49.00,11.50){$0$}
\put(110.00,11.50){$3$}
\put(150.00,50.00){$a_{i}b_{i}=b_{i}a_{i}=0,\qquad i=0,\ldots
,3,$} \put(150.00,35.00){$a_{i}a_{j}=b_{l}b_{m}$\quad for any\quad
$i,j,l,m\in\lbrace0,1,2,3\rbrace,$} \put(150.00,20.00){where
possible.} \put(59.50,57.00){$\circ$} \put(101.00,57.00){$\circ$}
\put(59.50,15.50){$\circ$} \put(101.00,15.50){$\circ$}
\put(69.00,62.00){\vector(1,0){27.00}}
\put(96.50,56.50){\vector(-1,0){27.00}}
\put(69.00,20.50){\vector(1,0){27.00}}
\put(96.50,15.00){\vector(-1,0){27.00}}
\put(106.50,52.00){\vector(0,-1){27.00}}
\put(64.50,25.00){\vector(0,1){27.00}}
\put(58.00,52.00){\vector(0,-1){27.00}}
\put(100.00,25.00){\vector(0,1){27.00}}
\end{picture}

\medskip

\noindent By Proposition~\ref{prop:simple}, this category has four
non-isomorphic simple modules $S_i$ where $i\in\{0,1,2,3\}$, with the
support in the given point $i$. The indecomposable modules were
described in \cite{bb:ir}. For the sake of completeness we repeat
here this classification.

We will treat the objects  of ${\bf B}$ as elements of $ \Z /4\Z$.
Consider the following three families of non-simple indecomposable
modules.

\medskip
\noindent {\it Finite family}.
Fix $i\in\{0,1,2,3\}$ and define the  ${\bf B}$-module  $M_i$
such that $M_i(j)=\Bbbk e_j$ for each  $j=0,1,2,3$ and
\ben
a_ie_i=e_{i+1}, \qquad
a_{i+1}e_{i+1}=e_{i+2},\qquad
b_{i-1}e_i=e_{i-2},     \qquad
b_{i-2}e_{i-1}=e_{i-2}
\end{equation*}
while $u_je_k=0$ for all other cases
of $u\in \{a,b\}$ and $j,k=0,\ldots ,3$.  Obviously,
$M_i$ is indecomposable module for any $i$.

\medskip
\noindent  {\it Infinite discrete families}.
Let $n\in \N$, $n>1$, and $j\in \Z_4$. Define
a ${\bf B}$-module $M_{n,j,1}$  (respectively, $M_{n,j,2}$) as follows.
Consider $n$ elements $e_1,\ldots ,e_{n}$.
A $\Bbbk$-basis of the vector space
$M_{n,j,1}(l)$ (respectively, $M_{n,j,2}(l)$) is the set of $e_k$
such that $j+k-1\equiv l(mod\,4).$
The elements $a_l$ and $b_{l-1}$ act as follows:

$$
a_{l}e_k =\left\{ \begin{array}{ll}
e_{k+1}, & \mbox{ if $l$ is even (resp., odd),  $k< n$ and
$j+k-1\equiv l(mod\,4)$;}\\
0, & \mbox{ otherwise.}
\end{array}
\right.
$$

$$
b_{l-1}e_k =\left\{ \begin{array}{ll}
e_{k-1}, & \mbox{ if $l$ is even (resp., odd),  $k>1$ and
$j+k-1\equiv l(mod\,4)$;}\\
0, & \mbox{ otherwise.}
\end{array}
\right.
$$
All modules $M_{n,j,1}$ and $M_{n,j,2}$, $n>1$, $0\leq j\leq 3$ are non-isomorphic
indecomposable ${\bf B}$-modules.

\medskip
\noindent  {\it Infinite continuous families}.
For each $\lambda \in \Bbbk$, $\la\ne 0$, and $d\in \Z$, $d>0$ define the
${\bf B}$-modules $M_{d,\lambda,1}$ and $M_{d,\lambda,2}$ as follows.
Set
\beql{mmod1}
\bal
M_{d,\lambda,1}(i)&=\Bbbk^{d},\\
M_{d,\lambda,1}(a_0)&=M_{d,\lambda,1}(a_2)=M_{d,\lambda,1}(b_1)={\bf I}_d,\\
M_{d,\lambda,1}(b_0)&=M_{d,\lambda,1}(b_2)=M_{d,\lambda,1}(a_1)=M_{d,\lambda,1}(a_3)=0, \\
M_{d,\lambda,1}(b_3)&=J_{d,\lambda}
\eal
\non
\end{equation}
and
\beql{mmod2}
\bal
M_{d, \lambda,2}(i)&=\Bbbk^{d}, \\
M_{d, \lambda,2}(b_0)&=M_{d, \lambda,2}(b_2)=M_{d, \lambda,2}(a_1)={\bf I}_d,\\
M_{d, \lambda,2}(a_0)&=M_{d, \lambda,2}(a_2)=
M_{d, \lambda,2}(b_1)=M_{d, \lambda,2}(b_3)=0, \\
M_{d, \lambda,2}(a_3)&=J_{d,\lambda},
\eal
\non
\end{equation}
where $J_{d,\lambda}$ is the Jordan cell of dimension $d$ with the eigenvalue $\lambda$.

All modules $M_{d, \lambda, k}$, $k=1,2$ are indecomposable and
corresponding indecomposable modules in $R^f_{\ell}$ have all
weight spaces of dimension $d$.

\bpr\label{prop:indecop} \bra\cite{bb:ir}, Proposition 3.3.1.\ket The
modules $S_i$,  $M_i$, $M_{n, i, 1}$, $M_{n, i, 2}$,
$M_{d,\lambda,1}$, $M_{d, \lambda,2}$ where $0\leq i\leq 3$, $d$
is a positive integer, $\lambda \in \Bbbk$, $\la\ne 0$, and $n\geq
2$ is an integer, constitute an exhaustive list of pairwise
non-isomorphic indecomposable ${\bf B}$-modules. \epr

The following theorem describes the representation type of
$R_{\ell}^f$.

\bth\label{thm:reptype}
\begin{itemize}
\item[(1)] If $|X_{\ell}|=0$ then $R^f_{\ell}$ is a
semisimple category with a unique indecomposable
{\rm(}=irreducible{\rm)} module;
\item[(2)] If
$|X_{\ell}|=1$ then $R^f_{\ell}$ has finite representation type;
\item[(3)] If $|X_{\ell}|=2$ then $R^f_{\ell}$ has
tame representation type if and only if $|I_k|=1$ for all $k\in
X$. Otherwise, $R^f_{\ell}$ has  wild representation type;
\item[(4)] If $|X_{\ell}|>2$ then $R^f_{\ell}$ has
wild representation type.
\end{itemize}
\eth

\begin{proof}
In the case when $|X_{\ell}|=1$ all indecomposable modules for
$S_{\ell}$ are described in Proposition~\ref{prop:indec}. Hence
$R^f_{\ell}$ has the finite representation type. If $|X_{\ell}|=2$ and
$|I_k|=1$ for each $k\in X$ then all indecomposable modules for
$S_{\ell}$ are described in Proposition~\ref{prop:indecop}. It
follows from the definition that $R^f_{\ell}$ has the tame
representation type in this case. If $|I_k|>1$ for at least one
$k$ then it is easy to construct a family of indecomposable
modules that depends on two continuous parameters. Hence, in this
case $R^f_{\ell}$ has the wild representation type. Suppose now that
$|X_{\ell}|>2$. Then $S_{\ell}$ contains a full subcategory of
wild representation type considered in \cite[Theorem~1]{bb:ir}. We
immediately conclude that $R^f_{\ell}$ has the wild representation
type. This completes the proof.
\end{proof}

\bco
\begin{enumerate}
\item If $|X_{\ell}|=0$ then the category
$R_{\ell}$ is a semisimple category with a unique indecomposable
module.
\item If $|X_{\ell}|=1$ then $R_{\ell}$
has finite representation type with indecomposable modules as in
Proposition~\ref{prop:indec}.
\end{enumerate}
\eco

\begin{proof}
Since cases $|X_{\ell}|\leq 1$ correspond  to the finite
representation type then the corresponding categories do not admit
infinite-dimensional indecomposable modules by \cite{a:rt} and
hence every indecomposable module belongs to $R^f_{\ell}$.
\end{proof}

\section*{Acknowledgment}

The first author is a Regular Associate of the ICTP and is
supported by the CNPq grant (Processo 300679/97-1). The first
author is grateful to the University of Sydney for support and hospitality.
 The second and
the third authors are grateful to FAPESP for the financial support
(Processos 2001/13973-0 and  2002/01866-7)  and to the University of S\~ao
Paulo for the warm hospitality during their visits.

\end{document}